\documentclass[12pt]{article}

\pdfoutput=1

\usepackage{amsmath}
\usepackage{graphicx}
\usepackage{enumerate}
\usepackage{natbib}
\usepackage{url} 

\newcommand{\blind}{1}

\usepackage{enumerate} 
\usepackage{cite}
\usepackage{authblk}
\usepackage[colorlinks,
linkcolor=blue,
anchorcolor=blue,
citecolor=blue]{hyperref}
\usepackage{natbib}
\setcitestyle{authoryear,round}

\usepackage{amssymb}
\usepackage{amsthm}
\usepackage{amsmath}
\usepackage{titlesec}
\usepackage{lineno,hyperref}
\modulolinenumbers[5]
\usepackage{stmaryrd}
\usepackage{dsfont}
\usepackage{enumerate}
\usepackage{graphicx,color}
\usepackage{float}
\usepackage{epstopdf}

\usepackage{subfigure}
\usepackage{amsfonts}
\usepackage{multirow}
\usepackage{bm}
\usepackage{appendix}
\usepackage{url}

\newtheorem{thm}{Theorem}

\newtheorem{remark}{Remark}[section]

\newtheorem{lem}{Lemma}[section]

\numberwithin{equation}{section}

\makeatletter  
\newif\if@restonecol  
\makeatother

\usepackage[ruled]{algorithm2e}
\usepackage[noend]{algpseudocode}

 \makeatother
\makeatletter
\@addtoreset{equation}{section}
\makeatother

\usepackage{xparse} 
\newcounter{example}[section]
\newenvironment{example}[1][]{\refstepcounter{example}\par\medskip
	\indent \textbf{Example~\theexample. #1} \rmfamily}{\medskip}

\usepackage{caption}
\begin{document}

	\def\spacingset#1{\renewcommand{\baselinestretch}%
		{#1}\small\normalsize} \spacingset{1}

	
	\if1\blind
	{
		\title{\bf Fast Calibration for Computer Models with Massive Physical Observations}
		\author{Shurui Lv\thanks{
				Correspondence should be addressed to Yan Wang (yanwang@bjut.edu.cn) 
			}\hspace{.2cm}\\
			School of Statistics and Data Science, Faculty of Science, Beijing University of Technology, Beijing, China\\
			Yan Wang \\
			School of Statistics and Data Science, Faculty of Science, Beijing University of Technology, Beijing, China\\
			and \\
			Jun Yu \\
			School of Mathematics and Statistics, Beijing Institute of Technology, Beijing 100081, China}
		\maketitle
	} \fi
	
	\if0\blind
	{
		\bigskip
		\bigskip
		\bigskip
		\begin{center}
			{\LARGE\bf Fast Calibration for Computer Models with Massive Physical Observations}
		\end{center}
		\medskip
	} \fi
	
	\bigskip
	\begin{abstract}
		Computer model calibration is a crucial step in building a reliable computer model. In the face of massive physical observations, a fast estimation for the calibration parameters is urgently needed. To alleviate the computational burden, we design a two-step algorithm to estimate the calibration parameters by employing the subsampling techniques. Compared with the current state-of-the-art calibration methods, the complexity of the proposed algorithm is greatly reduced without sacrificing too much accuracy. We prove the consistency and asymptotic normality of the proposed estimator. The form of the variance of the proposed estimation is also presented,  which provides a natural way to quantify the uncertainty of the calibration parameters. The obtained results of two numerical simulations and two real-case studies demonstrate the advantages of the proposed method.	 
		
		
	\end{abstract}
	
	\noindent%
	\vfill
	
	\newpage
	\spacingset{1.9} 
	
	\section{Introduction}
	\label{sec:intro}
	
	Computer models, or simulators, are being used increasingly to mimic complex systems in physics, engineering, and human processes. A computer model usually involves additional model parameters that cannot be determined or observed in physical processes. 
	Tuning these parameters is essential to match the computer outputs with the physical observations. 
	The corresponding stage is called the \emph{calibration} of computer models, and the parameters are usually referred as \emph{calibration parameters}.

	Several  methods have been established to estimate the calibration parameters, 
	such as the KO's calibration  \citep{kennedy2001bayesian}, $L_2$ calibration \citep{tuo2015efficient}, ordinary least squares (OLS) calibration \citep{wong2017frequentist}, orthogonal Gaussian Process (GP) model calibration \citep{plumlee2017bayesian}, projected kernel calibration \citep{tuo2019adjustments}, and Bayesian projected calibration \citep{xie2020bayesian}.  
	
	Despite the great success achieved by the aforementioned methods, the  computational complexity of those methods is greater than or equal to $O(n^2)$, where $n$ is the sample size of the physical observations.
	Consequently, there is an urgent need to develop a fast calibration method for dealing with the case that
	the available physical observations are massive, which is common in many calibration problems. 
	For example, in the calibration of a traffic flow model  \citep{hou2013calibration}, more than $30,000$ actual traffic data per year were collected by one detector in the UK's M25 freeway. All the detectors in the UK's M25 freeway can collect more than $10$ million data per year (\href{http://tris.highwaysengland.co.uk/detail/trafficflowdata}{http://tris.highwaysengland.co.uk/detail/trafficflowdata}). 
	
	Besides the common two computational barriers for general big data analysis: (1) data is too large to be held in a computer’s memory; and (2) computation by using whole physical observations is time-consuming in obtaining the estimated results.
	The calibration with massive physical observations meets additional challenge since that the computer models are always expensive. 
	To be precise, the computer outputs at all the physical observation points are needed to evaluate how well they match the real physical observations. 
	It is time-consuming to run computer models at all the physical observations. Although cheaper surrogate models can be built to mimic the computer outputs \citep{santner2018design}, it may still be unaffordable to predict the computer outputs at all the physical observations. 

	
	To address these computational barriers, the parallel computing method has been  successfully applied to the calibration for computer models, see  \citet{Cai2017model, chang2019computer,tsai2021calibration} and the references therein. 
	These methods improve computational efficiency by dividing the whole data set into subsets,  estimating the calibration parameters with the subsets at the same time, and combining the results from subsets to obtain a final estimator. 
	However, such a method is still expensive since it naturally requires more computational resources. Moreover, debugging is certainly easier and faster with fewer observations. 
	Besides parallel computing, another well-known tool to 
	reduce the computational burden is subsampling. The core of the subsampling method is to use non-uniform subsampling probabilities to include data points into the subsamples with a larger probability. It is worth mentioning that most of the current subsampling strategies focus on the specified regression models  \citep{mahoneystatistical, wang2019information, zhu2016gradient, wang2018optimal,yao2019optimal,wang2021optimal,ai2021optimal,yu2020optimal}, 
	{where the regression model 
		is a specified function on the parameters to be estimated.} 	 	
	However, in general, the computer models to be calibrated are  difficult to be represented as a specified function  of the calibration parameters. The existing subsampling methods can not be directly applied to the calibration for the computer models.
	
	{Nevertheless, the idea of subsampling is valuable to help researchers efficiently downsize the data volume and naturally accelerate the downstreaming analysis.}
	Similar to the existing subsampling literature, this work aims to select subsamples that contain more information about calibration parameters and derive a fast estimator based on the selected  subsamples.
	The contributions of this work include the following:
	\begin{itemize}
		\item We propose an inverse probability weighted least squares (IPWLS) calibration method using subsamples obtained based on the optimal subsampling probabilities.   A practical two-step algorithm is also given to approximate the IPWLS estimator with theoretical backups. 
		\item Compared to the existing lowest computational
		complexity method which is the OLS calibration, 
		the proposed calibration method is computational more efficient.
		\item The proposed estimator is proven to be consistent and asymptotically normal. Based on the asymptotic variance, a formula for the standard error of the proposed estimator is derived,  which provides a natural way for the uncertainty quantification of the calibration parameters.
	\end{itemize}
	
	The remainder of this work is organized as follows.  Section \ref{Section 2} gives a brief review of the OLS calibration for computer models, which has the lowest computational complexity among the existing calibration methods.  In Section \ref{Section 3}, we introduce an IPWLS calibration based on the OLS calibration. 
	We prove the asymptotic properties of the proposed estimator based on subsamples and derive the optimal subsampling probabilities based on A- and L-optimality criteria. In Section \ref{Section 4},  a practical two-step algorithm is proposed to 
	approximate the proposed estimator. The asymptotic results and the standard error of the resulting estimator are also presented.  
	Section \ref{Section 5} examines the proposed method through two numerical simulations and two real data studies. Section
	\ref{sec:conc} concludes the paper with some discussion.	All technique proofs are deferred in Supplementary Material.

	\section{Preliminaries}
	\label{Section 2}
	
	In this section, we give a brief introduction to  the OLS calibration \citep{wong2017frequentist}, 
	which enjoys the lowest computational complexity among all the existing calibration methods.

	Denote the input domain of physical experiments by $\bm\Omega$, which is assumed to be a convex and compact subset of $\mathbb R^d$. Let $\bm X=\{\bm{x}_1,\ldots,\bm x_{n}\}\subseteq\bm\Omega$ be the set of design points for the physical experiments and $\bm{Y}^p=(y^p_1,\ldots,y^p_n)^{T}$ be the corresponding physical responses, with the superscript $p$ indicating ``physical".  Suppose the physical experimental observation is generated by
	\begin{eqnarray}
		y^p_i=\zeta(\bm x_i)+e_i, \, i=1,\ldots,n,
		\label{eq2.1}
	\end{eqnarray}
	where $\zeta(\cdot)$  is an unknown function, called the \emph{true process}; and $e_i$'s are independent and identically distributed random variables with zero mean and finite variance $\sigma^2>0$. 
	
	Let $\bm\theta\in\bm\Theta$ to be the calibration parameter. Suppose that the parameter space $\bm\Theta$ is a compact subset of $\mathbb R^q$.
	Let $y^s(\bm x,\bm\theta)$  be the computer model which is used to simulate the true process $\zeta(\bm x)$. Here,  the superscript  $s$ indicates ``simulated".
	Following from \citet{tuo2015efficient} and \citet{wong2017frequentist}, the ``true" value of calibration parameter is defined as 
	\begin{eqnarray}
		\begin{aligned}
			\bm\theta^{*}:=\operatorname*{\arg\min}_{\theta \in \bm\Theta}\int_{\bm\Omega} \left[{\zeta(\bm x)-y^s(\bm x,\bm\theta)}\right]^2 d F(\bm x),
			\label{thetastar}
		\end{aligned}
	\end{eqnarray}
	where $F$ is the sampling distribution of $\bm x$.
	The goal of calibration is to estimate $\bm\theta^*$ such that the computer outputs are as close as possible to the physical experimental observations.
	A natural estimator for $\bm\theta^*$ is  the  minimizer of the following loss function  \citep{tuo2015efficient}
	\begin{eqnarray}
		\begin{aligned}
			\underset{\bm\theta \in \bm\Theta}{\arg\min} \, \frac{1}{n}\sum_{i=1}^{n}\left[y^{p}_i-y^{s}(\bm x_{i},\bm \theta)\right]^{2}. 
			\label{eq3}
		\end{aligned}
	\end{eqnarray}
	When the computer model is expensive,  only a limited number of computer experiments can be run. It is hard to perform computer experiments at all the physical observation points. Hence, a surrogate model is needed  to mimic the computer outputs.   
	The predictor of the computer model is denoted by $\hat y^s(\bm x,\bm \theta)$. By plugging $\hat y^s(\bm x,\bm \theta)$ into 	(\ref{eq3}),  the OLS estimator can be obtained as follows,
	\begin{eqnarray}
		\begin{aligned}
			&\hat{\bm \theta}:=\underset{\bm\theta \in \bm\Theta}{\arg\min} \, \frac{1}{n}\sum_{i=1}^{n}\left[y^{p}_i-\hat y^{s}(\bm x_{i},\bm \theta)\right]^{2}. 
			\label{eq3-expensive}
		\end{aligned}
	\end{eqnarray}
	
	By assuming the approximation error of the surrogate model $\hat y^s(\cdot,\cdot)$  is negligible, i.e., $\sup_{(\bm x,\bm\theta)\in \bm 
		\Omega\times \bm\Theta}|y^s(\bm x,\bm\theta)-\hat y^s (\bm x,\bm\theta)|=o_p(n^{-1/2})$, \citet{tuo2015efficient} and \citet{wong2017frequentist} proved that  the OLS estimator 	(\ref{eq3-expensive})  converges to the ``true" value of the calibration parameter	(\ref{thetastar})  with the convergence rate  $O_p(n^{-1/2})$ and 
	the OLS estimate  is asymptotically normal, i.e., $\hat{\bm\theta}\sim N(\bm\theta^*,\mathbf\Sigma)$ in distribution. Here, $\mathbf\Sigma$ can be evaluated by
	\begin{equation}
		\label{Sigma}
		\mathbf\Sigma=\frac{1}{n} {\mathbf{J}}^{-1}{\mathbf V}{\mathbf{J}}^{-1},
	\end{equation} 
	with
	\begin{equation}\label{J}
		\mathbf{J}=\int_{\bm\Omega}\frac{\partial^{2}[\zeta(\bm x)- y^{s}(\bm x, {\bm\theta}^*)]^{2}}{\partial\bm \theta \partial\bm \theta^{T}}d F(\bm x),
	\end{equation}
	and 
	\begin{equation}
		{\mathbf V}={4}\int_{\bm\Omega}\left\{[\zeta(\bm x)
		- y^{s}(\bm x,{\bm\theta}^*)]^{2}+\sigma^2\right\}\frac{\partial  y^{s}(\bm x,{\bm\theta}^*)}{\partial\bm\theta^{T}}\frac{\partial y^{s}(\bm x,{\bm\theta}^*)}{\partial\bm\theta}dF(\bm x).
	\end{equation} 
	
	From \citet{tuo2015efficient}, the OLS calibration is less efficient than  the $L_2$ calibration. Specifically,  the estimator given by the $L_2$ calibration method achieves the asymptotic variance with $\mathbf{\Sigma}_{opt}=\frac{1}{n}\mathbf{J}^{-1}\mathbf{V}_{opt}\mathbf {J}^{-1}$ and
	${\mathbf V}_{opt}={4}\sigma^2\int_{\bm\Omega}\frac{\partial  y^{s}(\bm x,{\bm\theta}^*)}{\partial\bm\theta^{T}}\frac{\partial y^{s}(\bm x,{\bm\theta}^*)}{\partial\bm\theta}dF(\bm x)$.
	Since the computer models are built under many assumptions and simplifications 
	the bias between  $\zeta(\cdot)$ and $y^s(\cdot,\bm\theta^*)$ can be large in practical situations \citep{kennedy2001bayesian}. 
	Suppose $\frac{\partial  y^{s}(\bm x,{\bm\theta}^*)}{\partial\bm\theta}\neq \bm 0$  for $\bm x\in\bm \Omega$, one can see that  $\mathbf{\Sigma}-\mathbf{\Sigma}_{opt}\geq 0$. This implies that the OLS calibration can lead to a larger variance compared with the $L_2$ calibration, especially for the scenario that the physical observations are limited. 
	Note that $\mathbf{\Sigma}-\mathbf{\Sigma}_{opt}=O_p(n^{-1})$.  The difference between $\mathbf{\Sigma}$ and $\mathbf{\Sigma}_{opt}$ is negligible as the physical observations increase.  For the calibration with massive physical observations,  the OLS calibration method is still considered as a method that can achieve the highest possible efficiency.
	
	Besides the high efficiency, the OLS estimator is conceptually clean and simple, easy to understand and calculate \citep{wong2017frequentist}.  For an expensive computer model,   suppose  a GP model \citep{santner2018design} is built as the surrogate by using $m$ computer experiments, the computational complexities of building the GP model and predicting the computer outputs at $\{(\bm x_i^T,\bm\theta^T)^T\}_{i=1}^n$ are $O(m^3)$ and $O(nm)$, respectively. Assume that the asymptotic step of finding the optimal solution of (\ref{eq3-expensive}) is $w_n$. For example, suppose  the gradient
	descent method is adopted to query the optimal solution of (\ref{eq3-expensive}), then $w_n=O(\log(n))$  \citep{bottou2010large}. 
	With each iteration, we need to recalculate the gradient of the empirical loss (\ref{eq3-expensive}) and re-predict the computer outputs at $\{(\bm x_i^T,\bm \theta^T)^T\}_{i=1}^n$. That is, the computational complexity of the OLS calibration is $O(w_n nq^2+w_n nm+m^3)$. 
	It is at least $O(n^3)$ lower than that of  the other existing calibration methods  because the computations of the other existing calibration methods involve the inverse of a covariance matrix of $n$ by $n$. 
	
	
	Although the OLS calibration method enjoys these advantages,  the computational problems are still the biggest snag when the physical observations are massive.  Even if we ignore the computing time to build the emulator for $y^s$ and to predict the computer outputs at the physical observation points, the computational complexity of the OLS calibration is still $O(w_n nq^2)$. Thus, a fast calibration method is needed to reduce the computational burden.

	\section{Calibration via Subsampling}
	\label{Section 3}
	To reduce the computational complexity,  we introduce a fast estimator by using 
	general Poisson subsampling method. 
	We prove that the proposed estimate can effectively approximate the true value of the calibration parameter 	(\ref{thetastar}).
	We also derive the optimal subsampling probabilities  in the sense that  the resulting estimate achieves the minimizing asymptotic mean square error (MSE).

	\subsection{IPWLS Calibration}
	
	Let $\pi_i$ to be the probability to sample the $i$-th data point for $i = 1,\ldots,n$. Let $S=\{(\bm x_{i}^{*T}, y^{p}(\bm x_{i}^{*}), \pi_{i} )^T \}_{i=1}^{r}$ be a set of subsamples 
	and the corresponding sampling  probabilities. To improve computational efficiency, a novel calibration method can be constructed by using the subsample set $S$  to approximate the full data OLS estimates. Inspired by sampling theory, we 
	suggest an  inverse
	probability weighted least squares estimator, denoted as \emph{IPWLS} estimate, by solving  the following minimization problem
	\begin{align}
		\tilde{\bm\theta}:=\underset{\bm\theta \in \bm\Theta}{\arg\min} \, l^{*}(\bm \theta), 
		\label{eq4}
	\end{align}	
	where
	\begin{align}
		l^{*}(\bm \theta)= \frac{1}{n}\sum_{i=1}^{r}\frac{1}{\pi_{i}}\left[y^{p}(\bm x_{i}^{*})-\hat y^{s}(\bm x_{i}^{*},\bm \theta)\right]^{2}.
		\label{eq5}
	\end{align}
	
	It can be seen that, with given subsample set $S$, the computer outputs need to be predicted only at the locations of the $r$ subsamples. Let the asymptotic number of iterations for optimizing (\ref{eq5}) is   $w_r$. For given sampling probabilities, the computational complexity of the proposed estimator is $O(w_rrq^2+w_rrm+m^3)$.
	By assuming $r\ll n$, which is natural in the big data setting, this method can drastically reduce the computational complexity.

	Now, we focus on the choice of the subsamples. As
	subsampling with replacement according to unequal probabilities requires accessing subsampling probabilities for the full data all at once. This takes a large memory to implement
	and may reduce the computational efficiency. To overcome this challenge, we apply the Poisson sampling \citep{yu2020optimal} to generate subsamples. 
	The general IPWLS calibration based on the    Poisson sampling is shown in Algorithm \ref{alg:poi}.
	\begin{algorithm}[htpb]
		\caption{General  IPWLS calibration}
		\label{alg:poi}   
		\SetAlgoNoLine
		\BlankLine
		\textbf{Initialization:} $S=\varnothing$;\\	
		\For{$i=1,...,n$} {  
			Generate a Bernoulli variable $a_{i} \sim B(1,\pi_{i})$;
			
			\uIf{$a_{i}=1$}  
			{  
				Update $S=S \cup  \left\{(\bm x_{i}, y_{i}^{p},\pi_{i} )\right\}$;
			}
			\uElse{
				$S=S$.
			}	
		}		
		\textbf{Estimation}: Obtain 
		$\tilde{\bm\theta}$ (\ref{eq4}) based on the subsample set $S$.
	\end{algorithm}  
	\begin{remark} 
		The subsample size, say $r^*$, in Algorithm \ref{alg:poi} is random such that $E(r^*)=\sum_{i=1}^{n}\pi_{i}$. We use $r=\sum_{i=1}^{n}\pi_{i}$ to denote the expected subsample size.
		As shown in \citet{AI2020101512} that $r$ is still concentrated around its expectation with a high probability, the sampling budget is still under control.
		
	\end{remark}

	From Algorithm \ref{alg:poi}, we can see that each $\pi_{i}$  can be calculated for each individual data point when scanning the full data. It can save computer memory effectively. Accordingly, the proposed calibration method addresses the three computational barriers faced by massive data calibration as mentioned in Section {\ref{sec:intro}}.

	Suppose the predictor $\hat y^s(\cdot,\cdot)$ is built by using the computer experimental data $\mathcal {D}^s= \left\{ (\bm x_i^s,\bm \theta_i^s);y_i^s \right\}_{i=1}^{m}$. 
	Following \citet{ezzat2018sequential}, we generate computer experimental design points independently of the physical observation points.  The accuracy of the predictor $\hat y^s(\cdot,\cdot)$ is independent on the physical observations and it is not required to be considered in the subsampling process. 

	\subsection{Asymptotic Properties}
	
	Denote the full physical data as $\mathcal{D}_n=\{(\bm x_i,y_i^p)\}_{i=1}^{n}$.  Assume the physical design points $\{\bm x_i\}_{i=1}^n$ are fixed and use $F_n$ to denote their empirical distribution.
	For two functions $f_1$ and $f_2$, let  $\|\cdot\|_n$ represents the $L_2(F_n) -$ norm, with  $<f_1,f_2>_n=\frac{1}{n}\sum_{i=1}^nf_1(\bm x_i)f_2(\bm x_i)$ and $\|\cdot\|$ represents the $L_2(F) -$ norm with $<f_1,f_2>=\int_{\bm\Omega}f_1(\bm x)f_2(\bm x)d F(\bm x)$.  Let  $\|\cdot\|_E$ to be the Euclidean norm, and  $\|f_1\|_{L_{\infty}(\bm\Omega)}=\sup_{\bm x\in\bm\Omega} f_1(\bm x)$. Next, we discuss the asymptotic properties of $\tilde{\bm\theta}$.  
	The approximating error of the predictor $\hat y^s(\cdot,\cdot)$  can be well-controlled by carefully selected design \citep{haaland2018framework} and  moderate amount of computer experiments. As suggested in \citet{tuo2015efficient} and \citet{wong2017frequentist}, we ignore the uncertainty of $\hat y^s(\cdot,\cdot)$ in this work by assuming that $\|\hat y^s-y^s\|_{L_{\infty}(\bm\Omega\times\bm\Theta)}=o_p(n^{-1/2})$. 
	
	Now, let us list the necessary assumptions for obtaining convergence results. 
	\begin{enumerate}[(H.1)]
		\item \label{a1}  $\{e_i\}$ is a sequence of i.i.d. random variables with zero mean and
		finite variance.  Also, assume that  ${\rm E}e_i^6 < \infty$.%
		\item \label{thetastar-in} $\bm\theta^*$
		is the unique solution to (\ref{thetastar}), and an interior point of $\bm\Theta$. 
		\item \label{designx} Suppose the design $\bm X$ satisfies that
		\begin{enumerate}[(i)]
			\item $\sup_{\bm\theta\in\bm\Theta} \left| \|\zeta-y^s(\cdot,\bm\theta)\|_n^2-\|\zeta-y^s(\cdot,\bm\theta)\|^2\right|=o_p(1)$;
			\item Elements of $\mathbf J_n-\mathbf J$ are $o_p(1)$, where  $\mathbf {J}$ is defined in (\ref{J}) and 
			\begin{equation}\label{Jn}
				\mathbf{J}_n=\frac{1}{n}\sum_{i=1}^n\frac{\partial^{2}[\zeta(\bm x_i)- y^{s}(\bm x_i,{\bm\theta}^*)]^{2}}{\partial\bm \theta \partial\bm \theta^{T}};
			\end{equation}
			\item $\frac{1}{n}\sum_{i=1}^n\frac{\partial[\zeta(\bm x_i)-y^{s}(\bm x_i,{\bm\theta}^*)]^{2}}{\partial\bm \theta}=O_p(n^{-1/2})$.		
		\end{enumerate}  
		\item \label{a4}  
		Assume $\lambda_{\min}({\mathbf J})>0$ and $\lambda_{\max}({\mathbf J})<\infty$; $\lambda_{\min}(\mathbf A)$ and $\lambda_{\max}(\mathbf A)$ represent the smallest and the largest eigenvalue of a matrix $\mathbf A$, respectively.
		\item \label{a2} Suppose that	{ $
			\left\|[ \zeta(\bm x)-y^{s}(\bm x,\bm\theta)]^6 \right\|_{L_{\infty}(\bm\Omega\times\bm\Theta)}< \infty$}.
		\item \label{a3} 	Define $\bm \Theta_0 \subset \bm\Theta$ as a neighborhood of ${\bm\theta}^*$.  Assume there are	
		\begin{enumerate}[(i)]
			\item $y^s(\bm x,\bm\theta)$ is three times continuously differentiable with respect to $\bm\theta$ in  $\bm\Theta_0$;  
			\item $ \frac{\partial  y^{s}(\bm x, {\bm\theta})}{\partial\theta_{j_1}}$ and $ \frac{\partial^{2}  y^{s}(\bm x, {\bm\theta})}{\partial\theta_{j_1} \partial\theta_{j_2}}$ are continuous with respect to $\bm x$ over  $\bm\Theta_0$;
			\item \label{a31}$\left\| \left[\frac{\partial  y^{s}(\bm x, {\bm\theta})}{\partial\theta_{j_1}} \right]^9\right\| _{L_{\infty}(\bm\Omega\times\bm\Theta_0)}< \infty$;
			\item \label{a32}$\left\|\left[\frac{\partial^{2}  y^{s}(\bm x, {\bm\theta})}{\partial\theta_{j_1} \partial\theta_{j_2}}\right]^{2}\right\| _{L_{\infty}(\bm\Omega\times\bm\Theta_0)} < \infty $; 
			\item \label{a33} $\left\|\left[\frac{\partial^{3} y^{s}(\bm x, {\bm\theta})}{\partial\theta_{j_1} \partial\theta_{j_2}\partial\theta_{j_3}}\right]^{2}\right\|_{L_{\infty}(\bm\Omega\times\bm\Theta_0)}< \infty$, $j_1,j_2,j_3=1,\ldots,q$.
		\end{enumerate} 
		\item \label{a6}Assume that $ y^{s}(\bm x_{i},\bm\theta)$ is  $m(\bm x_i)$-Lipschitz continuous. Exactly, $\forall$ $\bm\theta_1,\bm\theta_2 \in \bm\Theta$, there exist $m(\bm x_i)$ such that $
		|y^{s}(\bm x_{i},\bm\theta_1)- y^{s}(\bm x_{i},\bm\theta_2) | \leq m(\bm x_i)\left\|\bm\theta_1-\bm\theta_2 \right\|_E$. 
		\item \label{a5} $\underset{i=1,2,\ldots,n}{\max}(n\pi_{i})^{-1}=O_{p}(r^{-1})$.
	\end{enumerate}

	Assumption (H.\ref{a1}) bounds the observation errors. Assumption (H.\ref{thetastar-in}) ensures the uniqueness of the OLS estimate.
	Assumption (H.\ref{designx}) constrains the convergence of the physical observations. Assumption (H.\ref{a4}) ensures positive definiteness of the Hessian matrix. 
	Assumptions  (H.\ref{a2}) – (H.\ref{a6})  are constraints on the smoothness of the computer models. 
	Assumption (H.\ref{a5}) constrains the weights in (\ref{eq4}). It is mainly to avoid estimating equations being dominated by data points with extremely small sampling probability. It is common in classical sampling techniques \citep{breidt2000local}.  It can be seen that, relative to the conditions required for the convergence property of the OLS calibration \citep{wong2017frequentist}, only the constraint on the sampling probabilities is added here.

	\begin{thm}
		\label{th2:normal} 
		
		If Assumptions (H.\ref{a1}) – (H.\ref{a5}) hold, then as $n\rightarrow\infty$ and $r\rightarrow\infty$, there is, $\|\tilde{\bm\theta}-{\bm\theta}^*\|_E=o_p(r^{-1/2})$ and
		\begin{align}
			&\tilde{\bm\theta}-{\bm\theta}^*\to  N(\bm 0,\tilde{\mathbf {\Sigma}}^*) \label{eq14-2}
		\end{align}
		in distribution. Here, $\tilde{\mathbf {\Sigma}}^*=\mathbf\Sigma+{\rm E}_{\mathcal{D}_n }(\tilde {\mathbf \Sigma})$,  where $\mathbf\Sigma$ is the asymptotic variance-covariance matrix  of $\hat{\bm\theta}$, which is defined in (\ref{Sigma}) and  $\tilde {\mathbf \Sigma}=\tilde{\mathbf{J}}^{-1}\tilde{\mathbf V}\tilde{\mathbf{J}}^{-1}$ is the asymptotic variance-covariance matrix  of $\tilde{\bm\theta}$ conditional on $\mathcal{D}_n$, with 
		\begin{equation}
			\tilde{\mathbf{J}}=\frac{1}{n}\sum_{i=1}^n\frac{\partial^{2}[y^p_i- \hat y^{s}(\bm x_i,\hat{\bm\theta})]^{2}}{\partial\bm \theta \partial\bm \theta^{T}},
		\end{equation}
		and 
		\begin{equation}
			\tilde	{\mathbf V}=\frac{4}{n^{2}}\sum_{i=1}^{n}\frac{1-\pi_{i}}{\pi_{i}}[y^p_i
			-\hat y^{s}(\bm x_{i},\hat{\bm\theta})]^{2}\frac{\partial \hat y^{s}(\bm x_{i},\hat{\bm\theta})}{\partial\bm\theta^{T}}\frac{\partial  \hat y^{s}(\bm x_{i},\hat{\bm\theta})}{\partial\bm\theta}.
		\end{equation}
		
	\end{thm}
	Indeed, in the cases that the physical observation points do not satisfy the  moment conditions in Assumption (H.\ref{designx}), even though  $\tilde{\bm\theta}$ does not converge to the true value,   $\tilde{\bm\theta}$ can approximate the  full data estimate $\hat{\bm\theta}$ efficiently. 
	\begin{lem}
		\label{lemma1}
		If Assumptions (H.\ref{a1}) – (H.\ref{thetastar-in}) and Assumptions (H.\ref{a4}) – (H.\ref{a5}) hold, then as $n\rightarrow\infty$ and $r\rightarrow\infty$, there is, conditional on $\mathcal D_n$ in probability, 
		\begin{align}
			&\tilde {\mathbf \Sigma}^{-\frac{1}{2}}(\tilde{\bm\theta}-\hat{\bm\theta})\to  N(\bm 0,\mathbf I_q) 
		\end{align}
		in distribution. Here, $\mathbf I_q$ is the identity matrix of size $q$.
	\end{lem}

	From  Theorem \ref{th2:normal}, we know that  the uncertainty  of $\tilde{\bm\theta}$ can be divided into two parts. First is the uncertainty of $\hat{\bm\theta}$, which arises mainly from the physical observations and the misspecification of  the true process $\zeta(\cdot)$. The other one is the uncertainty of $\tilde{\bm\theta}$ given $\mathcal {D}_n$, which is from the subsampling.   Moreover,
	since $\|\hat{\bm\theta}-\bm\theta^*\|_E=O_p(n^{-1/2})$  and   $\|\tilde{\bm\theta}-\hat{\bm\theta}\|_E=O_p(r^{-1/2})$,  by assuming $r\ll n$, we have that  the uncertainty of $\tilde{\bm\theta}$ is mainly captured by the second part.    Hence, we focus the the uncertainty  arising from the subsampling and improve the accuracy of  $\tilde{\bm\theta}$ by choosing $\{\pi_i\}_{i=1}^n$ that makes some functional of $\tilde{\mathbf V}$ smaller.
	
	\subsection{Optimal Poisson Subsampling}
	As MSE is  a commonly used criterion for evaluating the performance of a parameter estimator, we query the 
	optimal subsampling probabilities  by minimizing the asymptotic mean square error (AMSE) of $\tilde{\bm\theta}$ given $\mathcal{D}_n$.  We abbreviate the criterion as mV-optimal. 
	It corresponds to the A-optimality criterion \citep{pukelsheim2006optimal} in the theory of optimal experimental design and is equivalent to minimizing the trace of $\tilde{\mathbf \Sigma}$ in Theorem \ref{th2:normal}. The following theorem gives the result.

	\begin{thm}
		\label{th3:pi} 
		In Algorithm \ref{alg:poi}, the subsampling strategy is mV-optimal if the subsampling probabilities are chosen as follows:
		\begin{equation}
			\pi_{i}^{mV}=r\frac{h_{i}^{mV}\wedge M}{\sum_{j=1}^{n} (h_{j}^{mV}\wedge M)}, \, i=1,\ldots,n,
			\label{pimV1}
		\end{equation}
		where 
		\begin{equation}
			h_{i}^{mV}=\left| y^p_i-\hat y^{s}(\bm x_{i},\hat{\bm\theta})\right|   \left[ \frac{\partial \hat y^{s}(\bm x_{i},\hat{\bm\theta})}{\partial\bm\theta}
			\tilde{\mathbf J}^{-2}\frac{\partial \hat y^{s}(\bm x_{i},\hat{\bm\theta})}{\partial\bm\theta^T} \right]^{\frac{1}{2}}.
		\end{equation}
		Let $h_{(1)}^{mV} \leq h_{(2)}^{mV} \leq \ldots \leq h_{(n)}^{mV}$ be the order statistics of $\{h_{i}^{mV}\}_{i=1}^{n}$. For convenience, denote $h_{(n+1)}^{mV}=\infty$ and assume that $h_{(n-r)}^{mV}>0$, then $M=\frac{1}{r-k}\sum_{i=1}^{n-k}h^{mV}_{(i)}$ and
		$k = \min \{ s \left.\right| 0\leq s\leq r, (r-s)h^{mV}_{(n-s)} < \sum_{i=1}^{n-s}h^{mV}_{(i)} \}$.
	\end{thm}

	As shown in Theorem \ref{th3:pi}, the optimal subsampling probabilities 
	$\{\pi_{i}^{mV}\}_{i=1}^{n}$ depend on data through both covariates and responses. 
	The inclusion of data points with lager values
	of $\left|y^p_i-\hat y^{s}(\bm x_{i},\hat{\bm\theta})\right|$ and $ \left[ \frac{\partial \hat y^{s}(\bm x_{i},\hat{\bm\theta})}{\partial\bm\theta}
	\tilde{\mathbf J}^{-2}\frac{\partial \hat y^{s}(\bm x_{i},\hat{\bm\theta})}{\partial\bm\theta^{T}} \right]$ will improve the robustness of the subsample estimator. A larger value of $\left|y^p_i-\hat y^{s}(\bm x_{i},\hat{\bm\theta})\right|$ means that the calibrated computer output differs more from the physical observation at $\bm x_i$.  Since $\|\hat{\bm\theta}-\bm\theta^*\|_E=o_p(n^{-1/2})$,  the distance between the physical observation and the calibrated computer output is dominated by the observation error and the model discrepancy $\zeta(\cdot)-y^s(\cdot,\bm\theta^*)$.
	The matrix $\frac{1}{n}\sum_{i=1}^n \left[ \frac{\partial \hat y^{s}(\bm x,\hat{\bm\theta})}{\partial\bm\theta^T}
	\frac{\partial \hat y^{s}(\bm x_{i},\hat{\bm\theta})}{\partial\bm\theta} \right] $ is the fisher information matrix of $\hat{\bm\theta}$  in the cases that the computer model $y^s(\cdot,\bm\theta^*)$ is perfect, i.e.,  $\zeta(\cdot)=y^s(\cdot,\bm\theta^*)$. 
	Therefore, intuitively, larger  $ \left[ \frac{\partial \hat y^{s}(\bm x_{i},\hat{\bm\theta})}{\partial\bm\theta}
	\tilde{\mathbf J}^{-2}\frac{\partial \hat y^{s}(\bm x_{i},\hat{\bm\theta})}{\partial\bm\theta^{T}} \right]$ indicates that $\bm x_i$  carries  more information about the calibration parameters.
	
	Evaluating the optimal subsampling probabilities requires the calculation of 
	$\tilde{\mathbf J}$, whose computational complexity is $O(q^2n)$. By adding the computational time of predicting the computer outputs and calculating the first derivative of the computer model at $\{(\bm x_i,\hat{\bm\theta})\}_{i=1}^n$,  the computational complexity of evaluating $\{h^{mV}_1,\ldots,h^{mV}_n\}$ is $O(q^2n+q^3n+m^3+mn)$ for given $\hat{\bm\theta}$.
	To save the computational time, following from \citet{wang2018optimal},  the optimal probabilities can be evaluated by minimizing $\rm tr(\tilde{\mathbf{V}})$.  We abbreviate this criterion as mVc-optimal. 
	It is called L-optimality in optimal experimental design \citep{pukelsheim2006optimal}.  The result is presented in the following theorem.
	
	\begin{thm}
		\label{th4:pi} 
		In Algorithm \ref{alg:poi}, the subsampling strategy is mVc-optimal if the subsampling probabilities are chosen as follows:
		\begin{equation}
			\pi_{i}^{mVc}=r\frac{h_{i}^{mVc}\wedge M}{\sum_{j=1}^{n} (h_{j}^{mVc}\wedge M)}, \, i=1,\ldots,n,
			\label{pi mVc}
		\end{equation}
		where 
		\begin{equation}
			h_{i}^{mVc}=\left| y^p_i-\hat y^{s}(\bm x_{i},\hat{\bm\theta})\right|   \left[ \frac{\partial \hat y^{s}(\bm x_{i},\hat{\bm\theta})}{\partial\bm\theta}
			\frac{\partial \hat y^{s}(\bm x_{i},\hat{\bm\theta})}{\partial\bm\theta^{T}} \right]^{\frac{1}{2}}.
		\end{equation}	
		Let $h_{(1)}^{mVc} \leq h_{(2)}^{mVc} \leq \ldots \leq h_{(n)}^{mVc}$ be the order statistics of $\{h_{i}^{mVc}\}_{i=1}^{n}$. For convenience, denote $h_{(n+1)}^{mVc}=\infty$ and assume that $h_{(n-r)}^{mVc}>0$, then
		$M=\frac{1}{r-k}\sum_{i=1}^{n-k}h^{mVc}_{(i)}$ 
		and $k = \min \{ s \left.\right| 0\leq s\leq r, (r-s)h^{mVc}_{(n-s)} < \sum_{i=1}^{n-s}h^{mVc}_{(i)} \}$.
	\end{thm}
	It can be seen that the alternative optimization criterion does  reduces the computing time since there is no need to calculate $\tilde{\mathbf{J}}$ and $\tilde{\mathbf{J}}^{-1}$.
	The computational complexity of evaluating $\{h^{mVc}_1,\ldots,h^{mVc}_n\}$ is reduced to $O(q^2n+m^3+mn)$ for given $\hat{\bm\theta}$.

	For ease of writing, we use a unified notation $\pi_i^{opt}$ to denote $\pi_i^{mV}$ and $\pi_i^{mVc}$ in Theorem \ref{th3:pi} and Theorem \ref{th4:pi}, that is,
	\begin{equation}
		\pi_{i}^{opt}=r\frac{h_{i}^{opt}\wedge M}{\sum_{j=1}^{n} (h_{j}^{opt}\wedge M)}=r\frac{h_{i}^{opt}\wedge M}{n\Psi}, \, i=1,\ldots,n,
		\label{pi opt1}
	\end{equation}
	where $M=\frac{1}{r-k}\sum_{i=1}^{n-k}h^{opt}_{(i)}$, $\Psi=\frac{1}{n}\sum_{j=1}^{n}(h_{j}^{opt}\wedge M)$, and 
	$h_{i}^{opt}$ is used to refer to $h_{i}^{mV}$ or $h_{i}^{mVc}$.  
	\begin{lem}
		\label{bestv}
		If the subsampling
		probabilities in Algorithm \ref{alg:poi} are chosen as $\pi_{i}^{opt}$, then  the  AMSE of  $\tilde{\bm\theta}$ given $\mathcal{D}_n$ can be represented by
		\begin{equation}
			\label{minsigma1}
			\frac{4}{n^2(r-k)}\left[\sum_{i=1}^{n-k} h_{(i)}^{opt}\right]\sum_{i=1}^{n-k} \left\{\frac{[h_{(i)}^{mV}]^2}{h_{(i)}^{opt} }\right\}-\frac{4}{n^2} \sum_{i=1}^{n-k}\left[h_{(i)}^{mV}\right]^2.
		\end{equation}
		Specifically, if the optimal subsampling
		probabilities is chosen as $\pi_{i}^{mV}$, then the 
		AMSE of  $\tilde{\bm\theta}$ given $\mathcal{D}_n$ is the minimum,  which can be represented by
		\begin{equation}
			\label{minsigma2}
			\frac{4}{n^2(r-k)}\left[\sum_{i=1}^{n-k} h_{(i)}^{mV}\right]^2-\frac{4}{n^2} \sum_{i=1}^{n-k}\left[h_{(i)}^{mV}\right]^2.
		\end{equation}
	\end{lem}
	
	As reported in \citet{yu2020optimal},  $k$ tends to be zero as $r/n$ becomes smaller. By setting $k=0$, we have that the minimum AMSE  of  $\tilde{\bm\theta}$
	depends on the mean and dispersion of $h_1^{mV},\ldots,h_n^{mV}$.  The full data with small and dispersed  $\{h_i^{mV}\}_{i=1}^n$  will deduce small minimum AMSE, if the subsamples are obtained according to the mV criterion.
	Since $\mathcal{D}_n$ is fixed before subsampling, this part is not discussed further. 
	Besides the full data, it can be seen that, the minimum AMSE depends on the subsamples size.
	To obtain more accurate parameter estimators, a larger amount of subsamples is  preferred.

	\section{Implementation of the IPWLS Calibration }
	\label{Section 4}
	
	The optimal probabilities in (\ref{pimV1}) and (\ref{pi mVc}) depend on $\hat{\bm\theta}$, $\tilde{\mathbf J}$, $M$ and $\Psi$, which are calculated by using the full data. Thus,
	an exact IPWLS estimator for the calibration parameters by using Algorithm \ref{alg:poi} can't be obtained directly. We propose a two-step algorithm to approximate
	the IPWLS estimator in this section. 
	
	\subsection{Two-step Algorithm}

	To practically implement the optimal probabilities, we need to replace $\hat{\bm\theta}$ and $\tilde{\mathbf J}$ with their respective pilot estimator $\tilde{\bm\theta}_0$ and $\tilde{\mathbf J}_0$, which can be obtained by a uniform subsample of $r_0$ with $r_0<n$. 
	Denote $S_{r_0}$ as the set of the pilot subsamples  and $|S_{r_0}|$ is the size of $S_{r_0}$. The pilot estimator $\tilde{\bm\theta}_0$ and $\tilde{\mathbf J}_0$ can be expressed respectively as
	\begin{equation}
		\label{theta0}
		\tilde{\bm\theta}_0=\underset{\bm\theta \in \bm\Theta}{\arg\min}\frac{1}{|S_{r_0}|}
		\sum_{S_{r_0}}[y_{i}^{p}-\hat y^{s}(\bm x_{i},{\bm\theta})]^2,
	\end{equation}
	and
	\begin{equation}
		\label{j0}
		\tilde{\mathbf J}_0=\frac{1}{|S_{r_0}|}\sum_{S_{r_0}}\frac{\partial^2[y_{i}^{p}-\hat y^{s}(\bm x_{i},\tilde{\bm\theta}_{0})]^2}{\partial\bm\theta\partial\bm\theta^{T}}.
	\end{equation}

	Furthermore,  to determine the subsampling probability of each data point separately, we use the pilot sample to approximate $M$ and $\Psi$. Following \citet{yu2020optimal}, choosing $M=\infty$ will not have a significant impact on the optimal subsampling probabilities as long as $r/n$ is small.  Thus, we set $M=\infty$ and  
	\begin{equation}
		{\Psi}_0=\frac{1}{|S_{r_0}|}\sum_{S_{r_0}}\left| y^p_i-\hat y^{s}(\bm x_{i},\tilde{\bm\theta}_0)\right| {\psi}_0(\bm x_i),
		\label{Psi}
	\end{equation}
	where   ${\psi}_0(\bm x_i)=\left[ \frac{\partial \hat y^{s}(\bm x_{i},\tilde{\bm\theta}_0)}{\partial\bm\theta}\tilde{\mathbf J}_0^{-2}\frac{\partial \hat y^{s}(\bm x_{i},\tilde{\bm\theta}_0)}{\partial\bm\theta^{T}} \right]^{1/2}$ for the mV criterion, and ${\psi}_0(\bm x_i)=\left[ \frac{\partial \hat y^{s}(\bm x_{i},\tilde{\bm\theta}_0)}{\partial\bm\theta}\frac{\partial \hat y^{s}(\bm x_{i},\tilde{\bm\theta}_0)}{\partial\bm\theta^{T}} \right]^{1/2}$ for the mVc criterion.
	
	Denote $\breve{\pi}_{i}^{opt}$ be the approximated subsampling probabilities with $\hat{\bm\theta}$, $\tilde{\mathbf J}$, $M$  $\Psi$, and $r$ in (\ref{pi opt1}) replaced by their  respective  pilot estimate, i.e. 
	\begin{equation}
		\breve{\pi}_{i}^{opt}=r\frac{\left| y_{i}^{p}-\hat y^{s}(\bm x_{i},\tilde{\bm\theta}_{0})\right|{\psi}_0(\bm x_i) }{n{\Psi}_0}, \, i=1,\ldots,n.
		\label{pi appro}
	\end{equation}
	
	{The estimator with $\breve{\pi}_{i}^{opt}$ inserted in (\ref{eq4}) may be sensitive to the data points with $y_{i}^{p}-\hat{y}^{s}(x_{i},\tilde{\bm\theta}_{0}) \approx 0$  if they are included in the subsamples. A small nugget can be added to $\breve{\pi}_{i}^{opt}$ to make the estimator more robust. Consider an extreme case, if $y_{i}^{p}-\hat{y}^{s}(x_{i},\tilde{\bm\theta}_{0})$ is zero for all $\bm x_i$, then $\tilde{\bm\theta}_{0}$ is accurate. In this case, the optimal subsampling probabilities should be  the  uniform subsampling probabilities that used to generate the pilot subsamples. As a result, we use the uniform subsampling probabilities as the nugget.
		Precisely, we use the following subsampling probabilities $\breve{\pi}_{i}^{w}$, which is a convex combination of $\breve{\pi}_{i}^{opt}$ in (\ref{pi appro}) and the uniform subsampling probabilities,}
	\begin{equation}
		\breve{\pi}_{i}^{w}=(1-\rho)r\frac{\left| y_{i}^{p}-\hat y^{s}(\bm x_{i},\tilde{\bm\theta}_{0})\right|\psi_0(\bm x_i)}{n{\Psi}_0}+\rho r\frac{1}{n}, \, i=1,\ldots,n,
		\label{pi weight}
	\end{equation}
	where $\rho\in(0,1)$. This weighted adjustment is also used in \citet{ma2014statistical} and \citet{yu2020optimal}. When $\rho$ is larger, the corresponding estimator will be more robust since $l^*(\bm\theta)$ in (\ref{eq4}) will not be inflated by
	data points with extremely small values of $\breve{\pi}_{i}^{opt}$.

	Considering some $\breve{\pi}_{i}^{w}$ may be larger than one due to approximated $\Psi$ and taking $M=\infty$, we need to use the inverse of $\breve{\pi}_{i}^{w}\wedge 1$'s as weights. 
	Since $\breve{\pi}_{i}^{w}\wedge 1$ depends only on  $\bm x_i$, $\bm y^p_i$, $\hat y^s(\bm x_i,\tilde{\bm\theta}_{0})$ and $\frac{\partial \hat y^{s}(\bm x_{i},\tilde{\bm\theta}_{0})}{\partial\bm\theta}$, each $\breve{\pi}_{i}^{w}$ can be calculated independently. That is, there is no need to calculate  all $\breve{\pi}_{i}^{w},i=1,\ldots, n$ at once.
	Thus, there is not need to load the full data into memory  and this is very computationally beneficial in terms of memory usage.
	The practical implementation above is summarized in Algorithm \ref{alg:imple}.
	\begin{algorithm}[!ht]
		\caption{Two-step Algorithm}
		\label{alg:imple}   
		\SetAlgoNoLine
		\BlankLine
		\textbf{Step 1:}
		Run Algorithm \ref{alg:poi} with subsample size $r_0$ and uniform subsampling probability $\bm\pi^{unif}=\{\frac{r_0}{n}\}_{i=1}^{n}$ to get a subsample set $S_{r_0}$. Use $S_{r_0}$ to obtain  the pilot estimates $\tilde{\bm\theta}_0$ from  (\ref{theta0}), $\tilde{\mathbf J}_0$ from (\ref{j0}),  and $\Psi_0$  from (\ref{Psi}).	
		
		\textbf{Step 2:} Set $S_0=S_{r_0}$ and $\pi_i=\breve{\pi}_i^w$, where $\breve{\pi}_i^w$ is defined in (\ref{pi weight}). Run Algorithm \ref{alg:poi}  to collect the subsample set $S_{r}=\left\{(\bm x^*_{i}, y^{p}(\bm x^*_{i}),\breve{\pi}_i^w \wedge 1 )\right\}_{i=1}^r$  and obtain
		$\breve{\bm\theta}$ based on the subsamples by solving  the following minimization problem  
		\begin{align}
			\breve{\bm\theta}:=\underset{\bm\theta \in \bm\Theta}{\arg\min} \, l_{\tilde{\bm\theta}_{0}}^{*}(\bm\theta), 
			\label{bre}
		\end{align}	
		where
		\begin{align}
			l_{\tilde{\bm\theta}_{0}}^{*}(\bm\theta)=\frac{1}{n}\sum_{i=1}^{r}\frac{1}{\breve{\pi}_{i}^{w}\wedge 1}\left[y^{p}(\bm x_{i}^{*})-\hat y^{s}(\bm x_{i}^{*},\bm \theta)\right]^{2}. 
			\label{breve theta}	
		\end{align}
	\end{algorithm}
	
	The optimal sampling probability  at each physical sample are required in Algorithm \ref{alg:imple}.  That is, the prediction  and the derivative of the computer outputs at $\{\bm x_i,\tilde{\bm\theta}_0\}_{i=1}^n$  is required  once in this Algorithm.  It leads that the  
	the computational complexity of evaluating $	\breve{\bm\theta}$ 
	becomes $O(nm+nq^3+w_r(r_0+r)q^2+w_r(r_0+r)m+m^3)$ for the mV criterion and $O(nm+nq^2+w_r(r_0+r)q^2+w_r(r_0+r)m+m^3)$ for the mVc criterion.
	
	\subsection{Asymptotic Properties}
	For the estimator obtained from Algorithm \ref{alg:imple}, 
	now we establish its consistency and asymptotic normality.
	
	\begin{thm}
		\label{th4:cvg} 
		If Assumptions (H.\ref{a1}) – (H.\ref{a5}) hold and $r_0r^{-\frac{1}{2}}\rightarrow 0$, then as $r\rightarrow \infty$ and $n\rightarrow \infty$, $\breve{\bm\theta}$ is consistent to ${\bm\theta}^*$
		in probability. 	
	\end{thm}

	In fact, as long as the pilot estimator $\tilde{\bm\theta}_0$ exists, Algorithm \ref{alg:imple} will produce a consistent estimator. We don't have to limit $r_0\rightarrow\infty$ in Theorem \ref{th4:cvg}. Of course, if $r_0\rightarrow\infty$, we can get from Theorem \ref{th2:normal} that $\tilde{\bm\theta}_0$ exists with probability approaching one.

	\begin{thm}
		\label{th6:normal} 
		If Assumptions (H.\ref{a1}) – (H.\ref{a5}) hold and $r_{0}r^{-\frac{1}{2}}\rightarrow 0$, then as $r_{0}\rightarrow \infty$, $r\rightarrow \infty$ and $n\rightarrow \infty$, there is $\|\breve{\bm\theta}-{\bm\theta}^*\|_E=o_p(r^{-1/2})$ and
		\begin{align}
			&\breve{\bm\theta}-{\bm\theta}^*\to  N(\bm 0,\breve{\mathbf {\Sigma}}^*) \label{eq15-2}
		\end{align}
		in distribution. Here, $\breve{\mathbf {\Sigma}}^*={\rm E}_{\mathcal{D}_n }(\breve {\mathbf \Sigma})+\mathbf\Sigma$,   where $\breve{\mathbf \Sigma }=\tilde{\mathbf J}^{-1}{\tilde{\mathbf V}^{w}}\tilde{\mathbf J}^{-1}$, and
		$${\tilde{\mathbf V}}^w=\frac{4}{n^{2}}\sum_{i=1}^{n}\frac{1-\pi_{i}^{w}\wedge 1}{\pi_{i}^{w}\wedge 1} [y^p_i-\hat y^{s}(\bm x_{i},\hat{\bm\theta})]^{2}\frac{\partial \hat y^{s}(\bm x_{i},\hat{\bm\theta})}{\partial\bm\theta^{T}}\frac{\partial \hat y^{s}(\bm x_{i},\hat{\bm\theta})}{\partial\bm\theta},$$ with 
		$$\pi_{i}^{w}=(1-\rho)r\frac{h_{i}^{mV}}{\sum_{j=1}^{n}h_{j}^{mV}}+\rho r\frac{1}{n},$$
		for the mV criterion and 
		$${\pi}_{i}^{w}=(1-\rho)r\frac{h_{i}^{mVc}}{\sum_{j=1}^{n}h_{j}^{mVc}}+\rho r\frac{1}{n},$$
		for the mVc criterion.	
	\end{thm}

	In Theorem 	\ref{th6:normal}, we require $r_{0}\rightarrow \infty$ to get a consistent pilot estimate which is
	used to identify the more informative data points in the second step. On the other hand,   $r_0$ should be much
	smaller than $r$ so that the more informative second step subsample dominates the weighted loss function. 

	\begin{lem}
		\label{lemma2}
		If Assumptions (H.\ref{a1}) – (H.\ref{thetastar-in}) and Assumptions (H.\ref{a4}) – (H.\ref{a5}) hold, and $r_{0}r^{-\frac{1}{2}}\rightarrow 0$, then 
		as $r_{0}\rightarrow \infty$, $r\rightarrow \infty$ and $n\rightarrow \infty$,  there is,				
		conditional on $\mathcal{D}_n$ in probability,
		\begin{align}
			&\breve{\mathbf \Sigma }^{-\frac{1}{2}}(\breve{\bm\theta}-\hat{\bm\theta}) \to  N(\bm 0,\mathbf I_q)
			\label{eq15}
		\end{align} 
		in distribution.
	\end{lem}
	
	\begin{lem} 
		Assume that $r\ll n$. The AMSE of  $\breve{\bm\theta}$ given $\mathcal{D}_n$ can be represented by
		\begin{equation}
			\label{brevesigmamv}
			\frac{4}{nr}\left[\frac{1}{n}\sum_{j=1}^n h_j^{mV}\right]\left[\sum_{i=1}^{n} \frac{(h_{i}^{mV})^2}{(1-\rho)h_{i}^{mV}+\rho\frac{1}{n}\sum_{j=1}^n h_j^{mV}}\right]-\frac{4}{n^2} \sum_{i=1}^{n}(h_{i}^{mV})^2
		\end{equation}
		for the mV criterion and 
		\begin{equation}
			\label{brevesigmamvc}
			\frac{4}{nr}\left[\frac{1}{n}\sum_{j=1}^n h_j^{mVc}\right]\left[\sum_{i=1}^{n} \frac{(h_{i}^{mV})^2}{(1-\rho)h_{i}^{mVc}+\rho\frac{1}{n}\sum_{j=1}^n h_j^{mVc}}\right]-\frac{4}{n^2} \sum_{i=1}^{n}(h_{i}^{mV})^2
		\end{equation}
		for the mVc criterion.
	\end{lem}

	By some easy calculations,  if  $\rho=1$, the 
	AMSE of  $\breve{\bm\theta}$  given $\mathcal{D}_n$ is $\frac{4}{r}\left[\frac{1}{n}\sum_{i=1}^n (h_i^{mV})^2\right]-\frac{4}{n^2} \sum_{i=1}^{n}(h_{i}^{mV})^2$.  If  $\rho=0$, the 
	AMSE of  $\breve{\bm\theta}$ given $\mathcal{D}_n$ is $\frac{4}{r}\left[\frac{1}{n}\sum_{i=1}^n h_i^{mV}\right]^2-\frac{4}{n^2} \sum_{i=1}^{n}(h_{i}^{mV})^2$ for the mV criterion.
	Since $\left[\frac{1}{n}\sum_{i=1}^n (h_i^{mV})^2\right]\geq \left[\frac{1}{n}\sum_{i=1}^n h_i^{mV}\right]^2$, together with Theorem 	\ref{th3:pi}, we have that, by carefully choosing $\rho$, the two-step algorithm is asymptotically more efficient than the uniform subsampling. That is,   $\tilde{\bm\theta}_0$ 
	is also asymptotically normal  from  Theorem \ref{th2:normal}, but  the trace of  its asymptotic variance is larger than that for the two-step  IPWLS estimator with the same subsample size.

	\subsection{Uncertainty Quantification of the IPWLS Estimator}
	\label{sec:error}
	Estimating the variance of the proposed estimator is crucial for statistical inferences such as hypothesis testing and confidence
	interval construction. As mentioned before, the uncertainty of the
	$\hat{\bm\theta}$ can be ignored when $r\ll n$. We use the variance of $\breve{\bm\theta}$ conditional on $\mathcal{D}_n$ to approximate its unconditional variance. When physical observations are limited, the bootstrap approach and Markov chain Monte Carlo (MCMC) method are widely used to estimate the variance of a frequency estimator \citep{wong2017frequentist} and a Bayesian estimator \citep{kennedy2001bayesian}, respectively.    Unfortunately,  the bootstrap method and MCMC method lose their advantages when the physical observations are massive due to their time-consuming computation. Thus, an explicit expression for the conditional variance of  $\breve{\bm\theta}$ is preferred to simplify the computation. 
	
	Note that, although we require $r\ll n$, in fact, the choice of $r$ can be relatively large. Thus we use the asymptotic variance-covariance matrix $\breve{\mathbf\Sigma}$ in Theorem \ref{th6:normal}  to get the estimated version.    This method is widely used in the uncertainty quantification for parameter estimation with large samples \citep{wang2018optimal}.
	This approach, however, requires calculations on the full data. 
	Following \citet{wang2018optimal},   variance-covariance matrix of $\breve{\bm\theta}$ given $\mathcal{D}_n$ can be estimated by using the subsamples $S_r$.  Specifically,  let
	${\Psi}_r=\frac{1}{r}\sum_{i=1}^r\left| y^p(\bm x_i^*)-\hat y^{s}(\bm x^*_{i},\tilde{\bm\theta}_0)\right| {\psi}_r(\bm x^*_i) 
	$ with ${\psi}_r(\bm x^*_i)=\left[ \frac{\partial \hat y^{s}(\bm x^*_{i},\tilde{\bm\theta}_0)}{\partial\bm\theta}\tilde{\mathbf J}_0^{-2}\frac{\partial \hat y^{s}(\bm x^*_{i},\tilde{\bm\theta}_0)}{\partial\bm\theta^{T}} \right]^{1/2}$ for the mV criterion, and ${\psi}_r(\bm x_i)=\left[ \frac{\partial \hat y^{s}(\bm x^*_{i},\tilde{\bm\theta}_0)}{\partial\bm\theta}\frac{\partial \hat y^{s}(\bm x^*_{i},\tilde{\bm\theta}_0)}{\partial\bm\theta^{T}} \right]^{1/2}$ for the mVc criterion.  We propose that $\breve{\mathbf \Sigma }$ can be approximated by 
	\begin{align}
		\breve{\mathbf \Sigma }_r=\breve{\mathbf J}_r^{-1}\breve{\mathbf V}_r^w\breve{\mathbf J}_r^{-1},
		\label{std err}
	\end{align}
	with
	\begin{equation}
		\breve{\mathbf{J}}_r=\frac{1}{n}\sum_{S_r}\frac{1}{{\breve\pi_{i}^{r}\wedge 1}}\frac{\partial^{2}[y^p(\bm x_i^*)-\hat y^{s}(\bm x^*_i,\breve{\bm\theta})]^{2}}{\partial\bm \theta \partial\bm \theta^{T}},
	\end{equation}
	and 
	\begin{equation}
		\breve	{\mathbf V}^w_r=\frac{4}{n^2}\sum_{S_r}\frac{1-{\breve\pi_{i}^{r}\wedge 1}}{({\breve\pi_{i}^{r}\wedge 1})^2}[y^p(\bm x_i^*)
		-\hat y^{s}(\bm x^*_{i},\breve{\bm\theta})]^{2}\frac{\partial \hat y^{s}(\bm x_i^*,\breve{\bm\theta})}{\partial\bm\theta^{T}}\frac{\partial \hat y^{s}(\bm x^*_{i},\breve{\bm\theta})}{\partial\bm\theta},
	\end{equation}
	where
	\begin{equation}
		\breve{\pi}_{i}^{r}=(1-\rho)r\frac{\left| y^{p}(\bm x_i^*)-\hat y^{s}(\bm x^*_{i},\tilde{\bm\theta}_0)\right|\psi_r(\bm x_i)}{n{\Psi}_r}+\rho r\frac{1}{n}.
		\label{pir}
	\end{equation}
	
	The above formula involves only the selected subsamples to estimate the conditional variance-covariance matrix. Based on the method of moments, 
	if $\breve{\bm\theta}$ is replaced by $\hat{\bm\theta}$, 
	$\breve	{\mathbf J}_r$ and $\breve	{\mathbf V}_r^w$
	are unbiased estimators of $\tilde	{\mathbf J}$ and $	{\tilde{\mathbf V}}^w$, respectively. 

	In addition to uncertainty quantification, the variance-covariance matrix   $\breve{\mathbf \Sigma }_r$ helps to determine the subsample size. Specifically, by
	pre-setting the width of the confidence band for parameter estimation, we can judge whether the current subsamples are sufficient. If the current confidence interval is too wide, the accuracy of parameter estimation can be achieved by adding more samples.

	\section{Numerical Studies}
	\label{Section 5}
	
	In this section, we evaluate the performance of the proposed estimator  through some simulations. Since the maximin  distance design is asymptotically optimum for building  $\hat y^s$ under a Bayesian setting \citep{johnson1990minimax},  we adopt the  Maximin Latin-hypercube design  \citep{santner2018design} to conduct the computer experiments. The common choice of $m$ is $m\geq 10(d+q)$, as recommended in \citet{Loeppky2009Special}.  To get a decent $\tilde{\bm\theta}_0$, following \citet{krishna2021robust}, we set $r_0$ on $2q+10d$.  To balance the robustness and accuracy of the proposed estimator,  following \citet{yu2020optimal,wang2018optimal}, we set $\rho=0.2$. The performance of the optimal subsampling probability is examined by 
	the following criterion:
	\begin{align*}
		{\rm RMSE}=\frac{1}{T}\sum_{t=1}^{T}\sum_{j=1}^q\left[\frac{\breve{\theta}^{(t)}_{j}-{\theta}_j^*}{\theta_j^*}\right]^2,
	\end{align*}
	where $\breve{\bm\theta}^{(t)}$ is the estimator obtained from the $t$-th subsamples with the optimal subsampling probabilities. We set $T=100$ throughout this section. 

	\subsection{Simulation Studies}
	
	\begin{example}
		Suppose the computer model is 
		\begin{align}
			y^s(x,\bm\theta) = 7[\sin(2\pi\theta_1-\pi)]^2 + 2(2\pi\theta_2-\pi)^2\sin(2\pi x -\pi),
			\label{model1}
		\end{align}
		where $x\in[0,1]$, $\bm\theta\in [0,0.25] \times [0,0.5]$.
		Suppose the true process $\zeta(x)=y^s(x,\bm\theta^*)$, with $\bm\theta^*=[0.2,0.3]^T$. The total number of samples $n$ is $10000$ and suppose the observation errors $e_i \sim N(0,0.2^2)$.
		
		Now we compare  the calculation time and accuracy of the proposed estimator  for three different subsampling criteria: uniform subsampling, mV and the mVc subsampling. The first step sample size $r_0$ is fixed on $14$, and  $r$ is from $\{100, 200, 300, 400, 600\}$. Since the uniform subsampling probability does not depend on unknown parameters and no pilot subsamples are required, it is implemented with subsample size $r_0+r$ for fair. 
		
		The calculation time for evaluating the estimator with different subsample sizes are shown in Table \ref{Table 1}. The calculation time for using the full data is also given for comparisons.
		\begin{table}[!ht]
			\centering
			\caption{Calculation time (seconds) v.s. different subsample sizes}
			\begin{tabular}{c|cccccc}
				\hline
				Criterion &  $r=100$&  $r=200$&  $r=300$&  $r=400$&  $r=600$&  $r=n=10000$\\
				\hline
				uniform&  $0.051$&  $0.076$&  $0.094$&  $0.117$&  $0.126$&  $0.952$\\
				mV&  $0.205$&  $0.227$&  $0.250$&  $0.270$&  $0.285$&  $1.379$\\
				mVc&  $0.188$&  $0.195$&  $0.204$&  $0.249$&  $0.257$&  $1.294$\\
				\hline
			\end{tabular}
			\label{Table 1}
		\end{table}
		It can be seen that by using the subsamples, the computational time to evaluate the estimator is greatly saved. The estimator based on the uniform subsampling method requires the least computing time since it does not need to calculate the subsampling probabilities. Because the optimal subsampling probabilities for the mV subsampling criterion require the calculation of 
		$\tilde{\mathbf J}$, evaluating the estimator based on the mV subsampling criterion is more time-consuming than that based on the mVc criterion.

		Figure \ref{E1 RMSE} shows the comparison of the accuracy for the proposed estimator.		 
		\begin{figure}[htbp]
			\centering
			\includegraphics[scale=0.4]{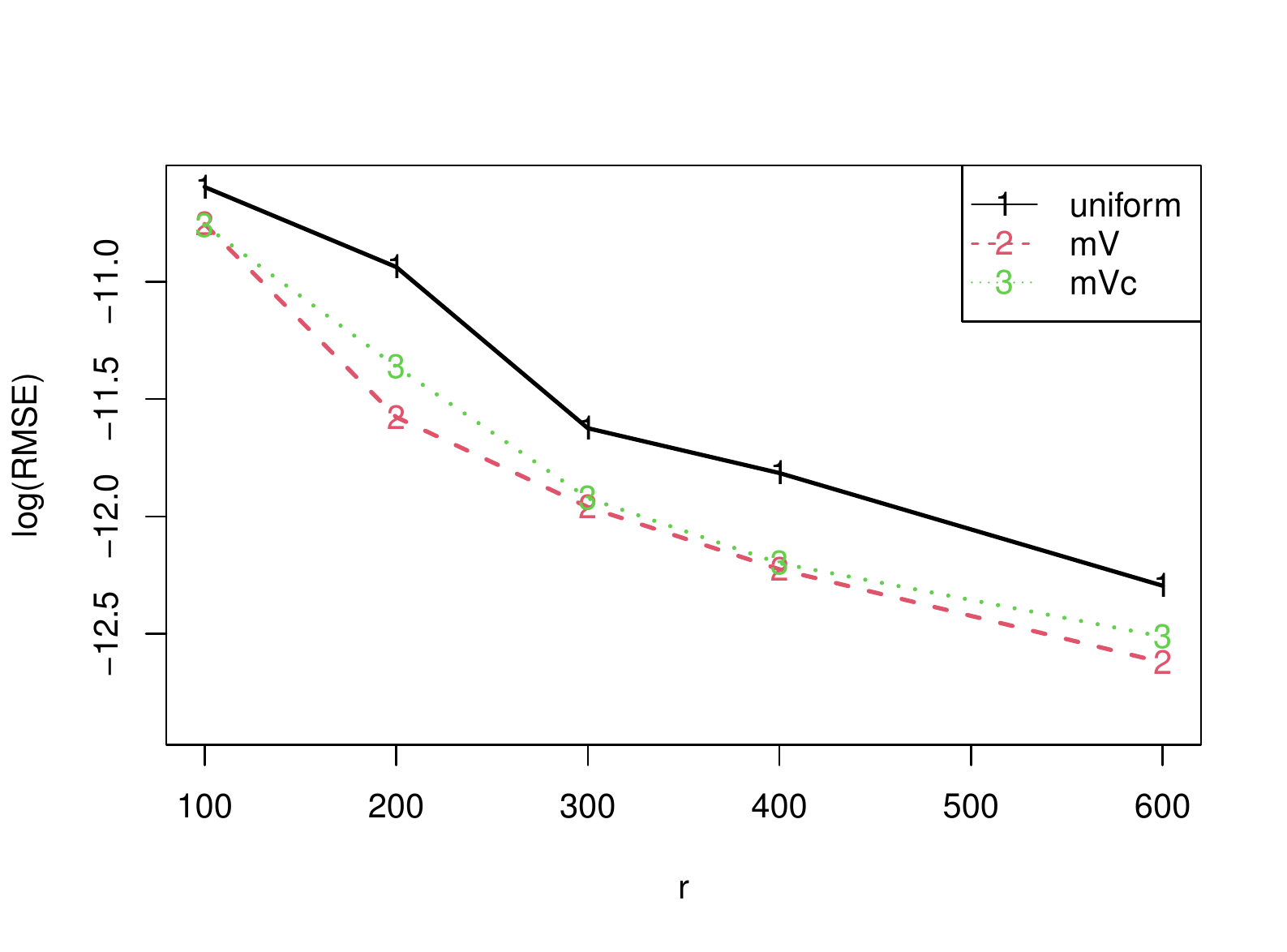} 
			\caption{Accuracy of the proposed estimator  for different subsample sizes $r$ and a fixed $r_0=14$ based on mV (red), mVc (green) and uniform (black) subsampling methods.}
			\label{E1 RMSE}
		\end{figure}
		The estimators by using the mV and mVc subsampling methods achieve smaller RMSEs. These results agree with the theoretical results in Section \ref{Section 3}. 
		
		The results of the relationship between subsamples and subsampling probabilities are presented in Figure \ref{x vs pi}. 
		\begin{figure}[!ht]
			\centering  
			\subfigure[$x_i^{mV}$ v.s. $\pi_i^{mV}$] {
				\begin{minipage}{0.4\textwidth}
					\centering        
					\includegraphics[scale=0.4]{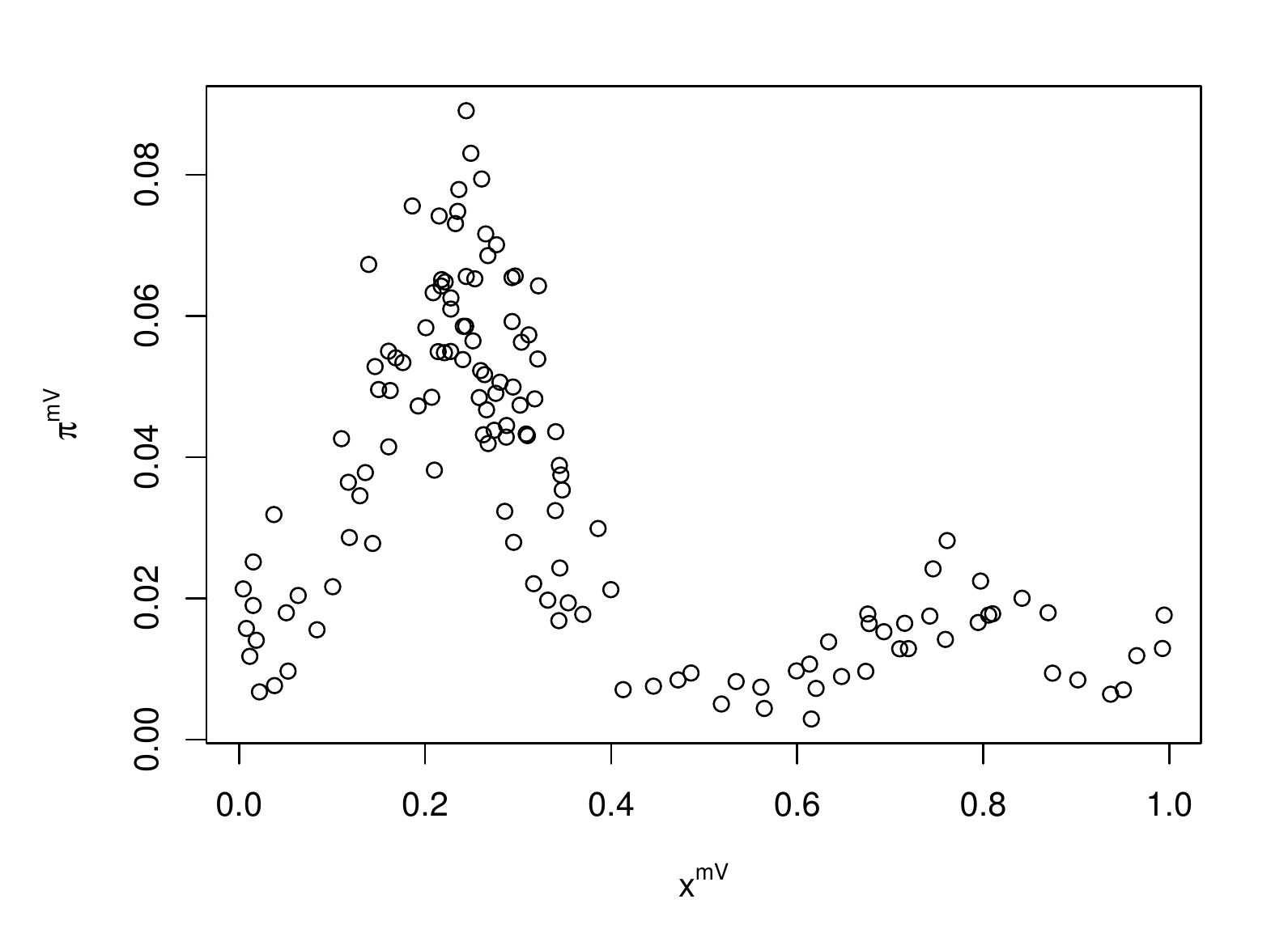}   
				\end{minipage}
			}
			\subfigure[$x_i^{mVc}$ v.s. $\pi_i^{mVc}$] {
				\begin{minipage}{0.4\textwidth}
					\centering      
					\includegraphics[scale=0.4]{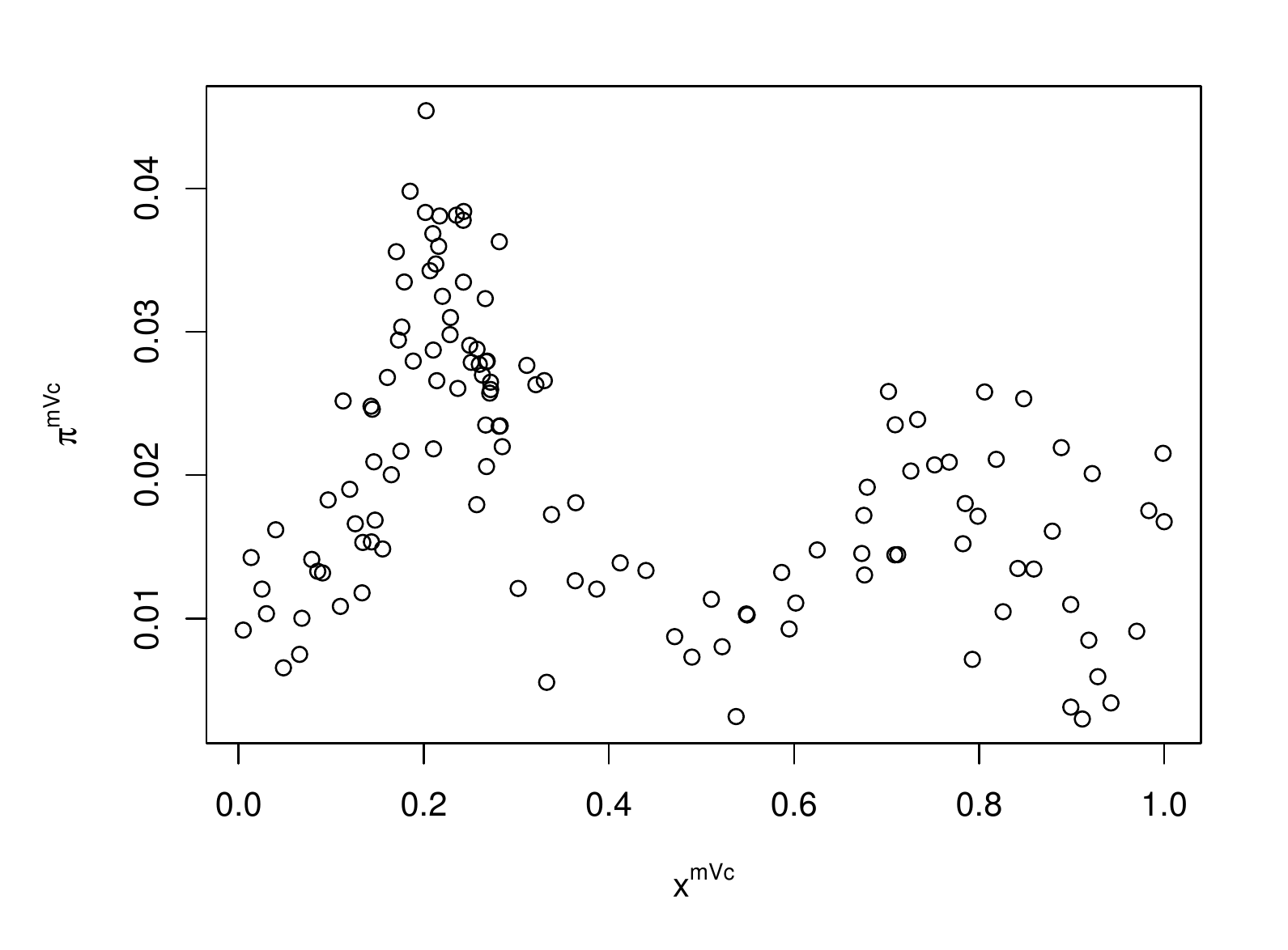}   
				\end{minipage}
			}
			\caption{The relationship between subsamples and subsampling probabilities based on mV and mVc methods with $r=100$.}
			\label{x vs pi}
		\end{figure}
		Combined with Figure \ref{E1 RMSE}, it states that the subsampling probabilities based on the mV and mVc criteria make full use of the information contained in the subsamples, which is not available in uniform subsampling, and this is why the two criteria have smaller RMSEs.
		In Figure \ref{x vs pi}, where $\pi_i^{mV}$ and $\pi_i^{mVc}$ are larger, the corresponding $\frac{\partial \hat y^{s}(\bm x_i^*,\breve{\bm\theta})}{\partial\bm\theta}$ is also larger, which is consistent with the results in Theorem \ref{th3:pi} and Theorem \ref{th4:pi}.

		To assess the performance of the proposed subsampling method for statistical inference,  the standard error given in (\ref{std err}) is used to estimate the variance-covariance matrices based on selected subsamples. We take $\theta_2$ as an example and present its coverage rates and lengths of $95\%$ confidence intervals in Table \ref{Table 5}. 
		\begin{table}[!ht]
			\centering
			\caption{Average lengths and coverage rates of $95\%$ confidence intervals for $\theta_2$.}
			\begin{tabular}{c|cccccc}
				\hline
				Criterion &  \multicolumn{2}{c}{uniform}&   \multicolumn{2}{c}{mV}&  \multicolumn{2}{c}{mVc}  \\
				&  length& coverage rate&   length& coverage rate&   length& coverage rate\\
				\hline
				$r=100$&  $0.0047$& $0.96$&  $0.0040$& $0.96$&  $0.0042$& $0.99$\\
				$r=200$&  $0.0034$& $0.96$&  $0.0029$& $0.98$&  $0.0029$& $0.94$\\
				$r=300$&  $0.0027$& $0.97$&  $0.0023$& $0.95$&  $0.0023$& $0.95$\\
				$r=400$&  $0.0024$& $0.98$&  $0.0019$& $0.94$&  $0.0020$& $0.96$\\
				$r=600$&  $0.0019$& $0.99$&  $0.0016$& $0.93$&  $0.0016$& $0.96$\\
				\hline
			\end{tabular}
			\label{Table 5}
		\end{table}
		Because the lengths of confidence intervals calculated by the uniform based subsampling method are the longest, its coverage rates at some sample sizes are slightly larger than that of the mV and mVc based subsampling methods.
		We can also obtain that the subsampling methods based on mV and mVc require fewer subsamples than that based on uniform to achieve considerable statistical inference effect. In general, mV and mVc subsampling methods outperform the uniform subsampling method. 
	\end{example}

	\begin{example}
		Suppose $\zeta(\cdot)$ can be expresses as 
		\begin{equation}
			\zeta(\bm x) = \frac{x_1}{2}\left[ \sqrt{1+(x_1+x_3^2)\frac{x_4}{x_1^2}}-1 \right] + (x_1+3x_4)\exp[1+\sin(x_3)]
			\notag
		\end{equation}
		and the computer model \citep{xiong2013sequential} is
		\begin{equation}
			y^s(\bm x,\bm \theta) = (\theta_1+\frac{\sin(x_1)}{10})\zeta(\bm x) + \theta_2(-2x_1+x_2^2+x_3^2)+0.5,
			\notag
		\end{equation}
		where $\theta_1$ and $\theta_2$ are two calibration parameters, with $\bm\theta \in [-5,5]^2$, and $\bm x\in[0,1]^4$. Let $\bm X=\{\bm x_1,\bm x_2,\ldots,\bm x_n\}$ be the set of design points in a MmLHDs, with $n=10000$; $e_i$'s are mutually independent and follow $N(0,0.1^2)$. We apply BB algorithm \citep{varadhan2010bb} to find that the local optimal point $\bm\theta^{*} = [0.895, 0.267]^T$. 
		To investigate the performance of the proposed method,  we select the subsamples with the subsample size equals to 400, 500, 600, 700, and 800, respectively.

		The computational time for different subsample sizes under the three subsampling criteria are summarized in Table \ref{Table 2}.
		\begin{table}[!ht]
			\centering
			\caption{Calculation time (seconds) v.s. different subsample sizes}
			\begin{tabular}{c|cccccc}
				\hline
				Criterion &  $r=400$&  $r=500$&  $r=600$&  $r=800$&  $r=900$&  $r=n=10000$\\
				\hline
				uniform&  $0.150$&  $0.153$&  $0.168$&  $0.179$&  $0.237$&  $1.517$\\
				mV&  $0.360$&  $0.371$&  $0.429$&  $0.492$&  $0.5338$&  $1.581$\\
				mVc&  $0.316$&  $0.360$&  $0.367$&  $0.458$&  $0.493$&  $1.531$\\
				\hline
			\end{tabular}
			\label{Table 2}
		\end{table}
		To assess the performance of Algorithm \ref{alg:imple} under three subsampling probabilities,  we also evaluate the $\rm RMSE$ for the three different subsampling criteria.
		Figure \ref{E2 RMSE} compares the accuracy of different estimators.
		\begin{figure}[htbp]
			\centering
			\includegraphics[scale=0.4]{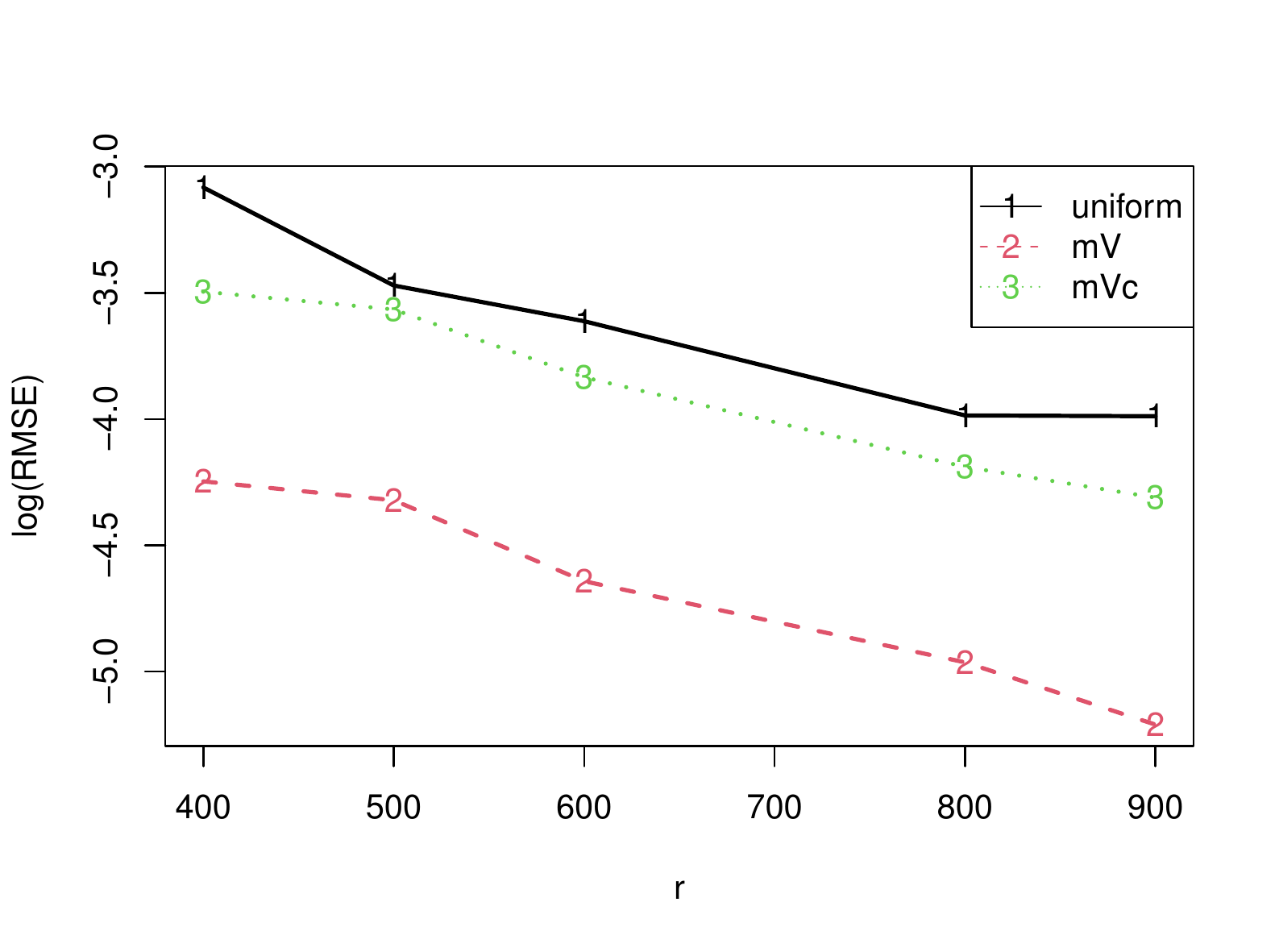} \label{3}
			\caption{Accuracy of the proposed estimator  for different subsample sizes $r$ and a fixed $r_0=44$ based on mV (red), mVc (green) and uniform (black) subsampling methods.}
			\label{E2 RMSE}
		\end{figure}
		It can be concluded that the estimators obtained by using the mV and mVc subsampling methods outperform that obtained by using the uniform subsampling method. Otherwise, RMSE  decreases as $r$ increases, which confirms the theoretical result on the consistency of the subsampling methods.
		
		Combining Figure \ref{E2 RMSE} and Figure \ref{i vs pi}, we can also conclude that the uniform subsampling method does not make full use of the information contained in the subsamples, so it is not as good as the subsampling method based on mV and mVc.
		\begin{figure}[!ht]
			\centering 
			\subfigure[$i$ v.s. $\pi_i^{mV}$] {
				\begin{minipage}{0.4\textwidth}
					\centering        
					\includegraphics[scale=0.4]{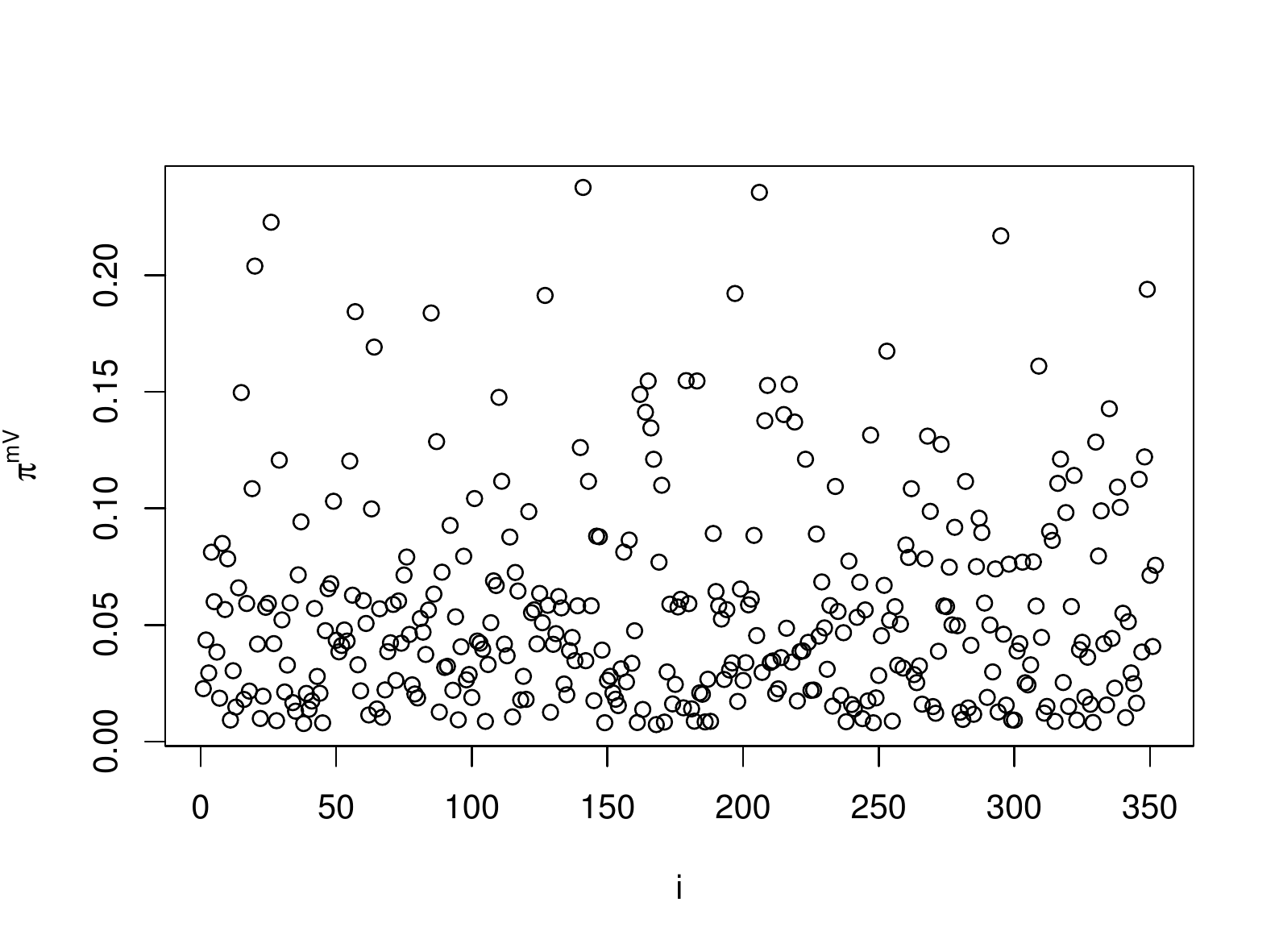}   
				\end{minipage}
			}
			\subfigure[$i$ v.s. $\pi_i^{mVc}$] {
				\begin{minipage}{0.4\textwidth}
					\centering      
					\includegraphics[scale=0.4]{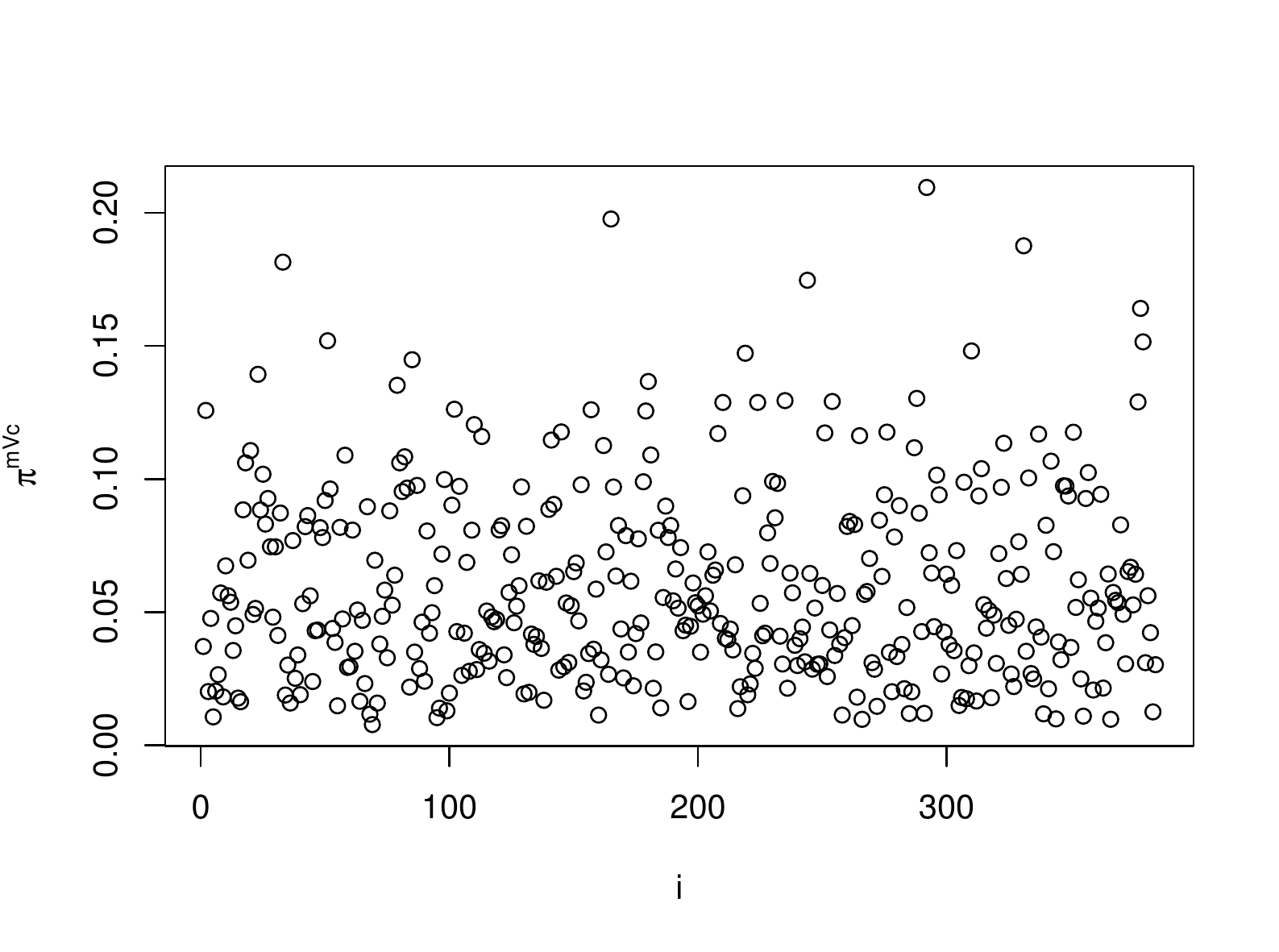}   
				\end{minipage}
			}
			\caption{The distribution of subsampling probabilities based on mV and mVc methods with $r=400$.}
			\label{i vs pi}
		\end{figure}

		Similarly, we get $95\%$ confidence intervals of parameters based on uniform, mV and mVc subsampling methods by simple calculation. As an example, we take $\theta_2$ as the parameter of interest. The average lengths and coverage rates of confidence intervals are reported in Table \ref{Table 6}. 
		\begin{table}[!ht]
			\centering
			\caption{Average lengths and coverage rates of $95\%$ confidence intervals for $\theta_2$.}
			\begin{tabular}{c|cccccc}
				\hline
				Criterion &  \multicolumn{2}{c}{uniform}&   \multicolumn{2}{c}{mV}&  \multicolumn{2}{c}{mVc}  \\
				&  length& coverage rate&   length& coverage rate&   length& coverage rate\\
				\hline
				$r=400$&  $0.208$& $0.92$&  $0.135$& $0.96$&  $0.180$& $0.93$\\
				$r=500$&  $0.187$& $0.97$&  $0.117$& $0.96$&  $0.166$& $0.95$\\
				$r=600$&  $0.169$& $0.96$&  $0.107$& $0.97$&  $0.149$& $0.93$\\
				$r=800$&  $0.148$& $0.95$&  $0.090$& $0.95$&  $0.129$& $0.94$\\
				$r=900$&  $0.139$& $0.96$&  $0.083$& $0.97$&  $0.119$& $0.95$\\
				\hline
			\end{tabular}
			\label{Table 6}
		\end{table}
		It shows that mV and mVc based subsampling methods have similar performances and they both perform better than uniform subsampling method. At some sample points, the coverage rates of uniform based subsampling method is greater than that of mVc based subsampling method because its corresponding confidence interval length is the longest.  As $r$ increases, length of confidence interval decreases. All of these confirm that the proposed asymptotic variance-covariance formula in (\ref{std err}) works well. 
	\end{example}

	\subsection{Real Case Studies}
	
	\begin{example}
		Well-calibrated traffic flow models are of great help in solving road congestion problems. We apply the proposed method to calibrate  the following dual-regime modified Greenshields traffic flow model \citep{hou2013calibration} 
		\begin{equation}
			\label{greenshields}
			v_{i}=\left\{\begin{array}{ll}
				u_{f}, & 0<k_{i}<k_{\mathrm{bp}}, \\
				v_{0}+\left(v_{f}-v_{0}\right)\left(1-\frac{k_{i}}{k_{\mathrm{jam}}}\right)^{\alpha}, & k_{\mathrm{bp}}<k_{i}<k_{\mathrm{jam}},
			\end{array}\right.
		\end{equation}
		where $v_i$ is the vehicle speed, and $k_i$ is the total carriageway density. A total of six parameters affect the shape of the model, which are breakpoint density $k_{bp} 
		\in [0,20]$, free-flow speed $u_f \in [100,120]$, 
		speed intercept $v_f \in [150,220]$, minimum speed $v_0 \in [0,10]$, jam density $k_{jam}\in [200,250]$ and shape parameter $\alpha \in [0,10]$. 
		This model is widely used to predict and explain the trends that are observed in real freeways traffic flows. The traffic data are collected  by the detector at London Orbital Motorway M25/4883A on link 199131002, located near Heathrow Airport.
		We select $35040$ traffic data to estimate the calibration parameters. Details about the data can be found in \href{http://tris.highwaysengland.co.uk/detail/trafficflowdata}{http://tris.highwaysengland.co.uk/detail/trafficflowdata}.
		
		Likewise, comparison of the calculation time for different subsample sizes is shown in Table \ref{Table 3}.
		\begin{table}[!ht]
			\centering
			\caption{Calculation time (seconds) v.s. different subsample sizes}
			\begin{tabular}{c|cccccc}
				\hline
				Criterion &  $r=100$&  $r=200$&  $r=400$&  $r=600$&  $r=800$&  $r=n=35040$\\
				\hline
				uniform&  $0.153$&  $0.159$&  $0.217$&  $0.321$&  $0.353$&  $3.227$\\
				mV&  $0.476$&  $0.602$&  $0.684$&  $0.716$&  $0.732$&  $3.583$\\
				mVc&  $0.443$&  $0.549$&  $0.573$&  $0.622$&  $0.689$&  $3.352$\\
				\hline
			\end{tabular}
			\label{Table 3}
		\end{table}
		Since the true values of the calibration parameters are unknown,  the difference between the subsamples estimator and the full data estimator is used to compare the performance of different estimators. Accordingly, define the following,
		\begin{align}
			\label{RMSE-HAT}
			{\rm RMSE_f}=\frac{1}{T}\sum_{t=1}^{T}\sum_{j=1}^q\left[\frac{\breve{\theta}^{(t)}_{j}-\hat{\theta}_j}{\hat{\theta}_j}\right]^2.
		\end{align}
		The relationship between  $\rm RMES_f$ and $r$ is presented in Figure \ref{MSE3}.
		\begin{figure}[!ht]
			\centering
			\includegraphics[scale=0.4]{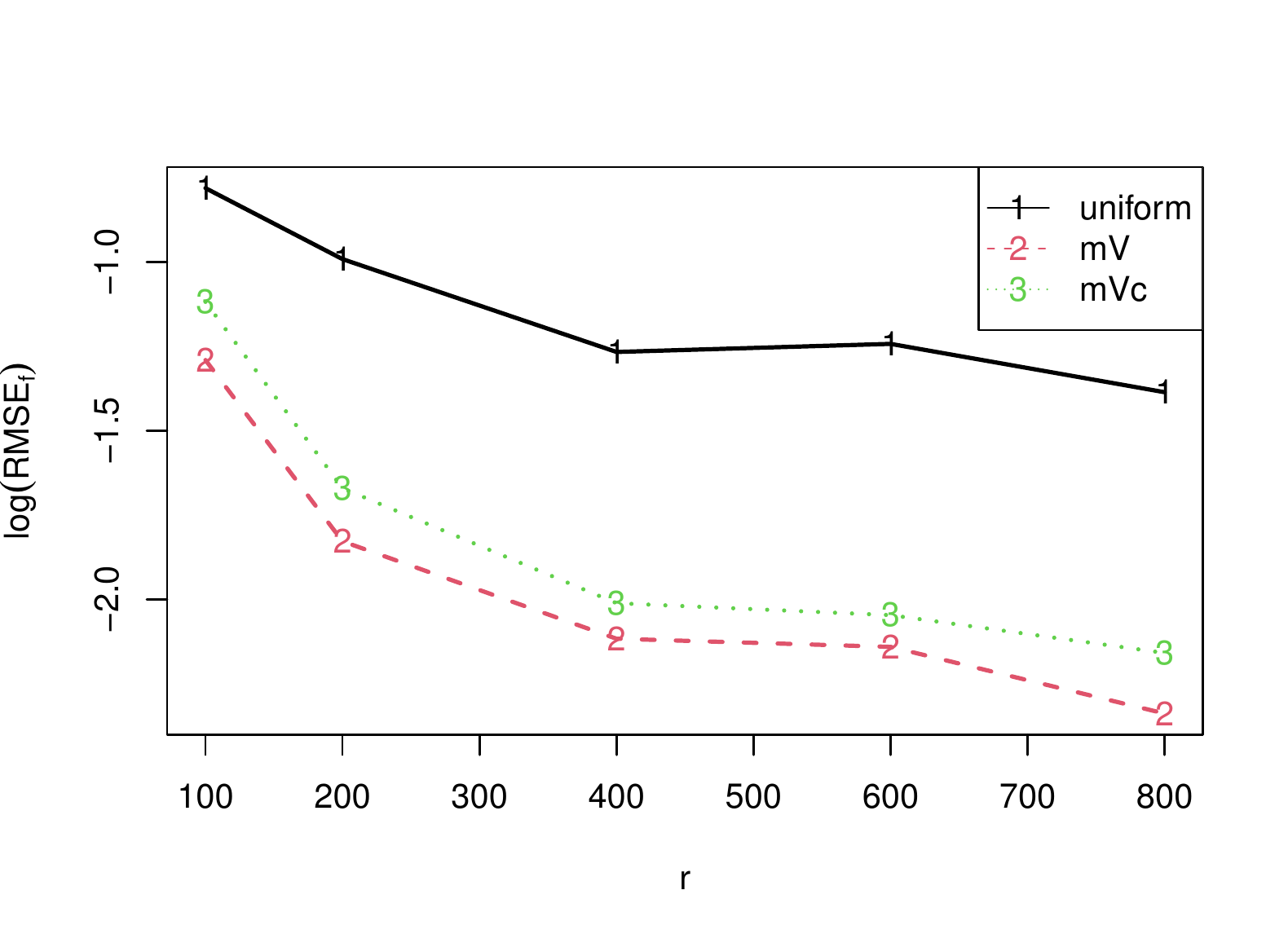}
			\caption{Accuracy of the proposed estimator  for different subsample sizes $r$ and a fixed $r_0=22$ based on mV (red), mVc (green) and uniform (black) subsampling methods.}
			\label{MSE3}
		\end{figure}
		It shows that the estimators  based on the mV and mVc criteria outperform the estimator based on the uniform subsampling method. 
		
		Next we take $\alpha$ as an example and present its $95\%$ confidence intervals by different subsampling methods. The average lengths of the $95\%$ confidence intervals are presented in Table \ref{Table 7}. We do not have a true value in the actual data. As a consequence, the coverage rates are not reported here. So is the next example.
		\begin{table}[!ht]
			\centering
			\caption{Average lengths of $95\%$ confidence intervals for $\alpha$.}
			\begin{tabular}{c|cccccc}
				\hline
				Criterion &  $r=100$&  $r=200$&  $r=400$&  $r=600$&  $r=800$\\
				\hline
				uniform&  $11.364$&  $7.761$&  $5.259$&  $4.561$&  $3.950$\\
				mV&  $3.623$&  $2.460$&  $2.263$&  $2.250$&  $1.920$\\
				mVc&  $3.765$&  $2.687$&  $2.106$&  $1.894$&  $1.470$\\
				\hline
			\end{tabular}
			\label{Table 7}
		\end{table}
		For some sample sizes, the confidence interval lengths obtained by the mV subsampling method are longer than that obtained by the mVc subsampling method. It is acceptable as
		the optimal subsampling probabilities based on mV is obtained by minimizing $\rm tr(\mathbf V)$,
		that is, minimizing the sum of variance of parameters. 
		In summary, better statistical inference effect can be achieved by the mV and mVc based subsampling methods rather than the uniform based subsampling methods.
	\end{example}

	\begin{example}
		In this example, we consider the calibration of the expensive  VarKarst-R model \citep{hartmann2015large}, which is a large-scale simulation model to assess karstic groundwater recharge.
		The model is a  function of precipitation and potential evapotranspiration and its output is karst recharge and actual evapotranspiration. 
		It contains four calibration parameters, which are variability constant $a \in [3,5]$, mean epikarst storage coefficient $K_{epi} \in [25,30]$, mean soil storage capacity $V_{soil} \in [530,545]$, and mean epikarst storage capacity $V_{epi} \in [400,430]$. A well-calibrated model can help to improve karst water budgets, which inform decisions on drinking water supply and flood risk management.
		We select $1488$ physical observations from December 20, 2003 to January 15, 2008, which is available at \href{https://github.com/KarstHub/VarKarst-R-2015}{https://github.com/KarstHub/VarKarst-R-2015}. 
		
		Comparison of the calculation time for different subsample sizes is shown in Table \ref{Table 4}. 
		\begin{table}[!ht]
			\centering
			\caption{Calculation time (seconds) v.s. different subsample sizes}
			\begin{tabular}{c|cccccc}
				\hline
				Criterion &  $r=300$&  $r=400$&  $r=500$&  $r=600$&  $r=800$&  $r=n=1488$\\
				\hline
				uniform&  $30.452$&  $43.346$&  $47.677$&  $48.108$&  $64.320$&  $89.244$\\
				mV&  $31.776$&  $48.869$&  $49.411$&  $56.166$&  $66.198$&  $94.782$\\
				mVc&  $31.143$&  $46.470$&  $48.267$&  $53.722$&  $65.862$&  $92.562$\\
				\hline
			\end{tabular}
			\label{Table 4}
		\end{table}
		The $\rm RMSE_f$s with a fixed $r_0=28$ and different subsample sizes are shown in Figure \ref{MSE4}. 
		\begin{figure}[htbp]
			\centering
			\includegraphics[scale=0.4]{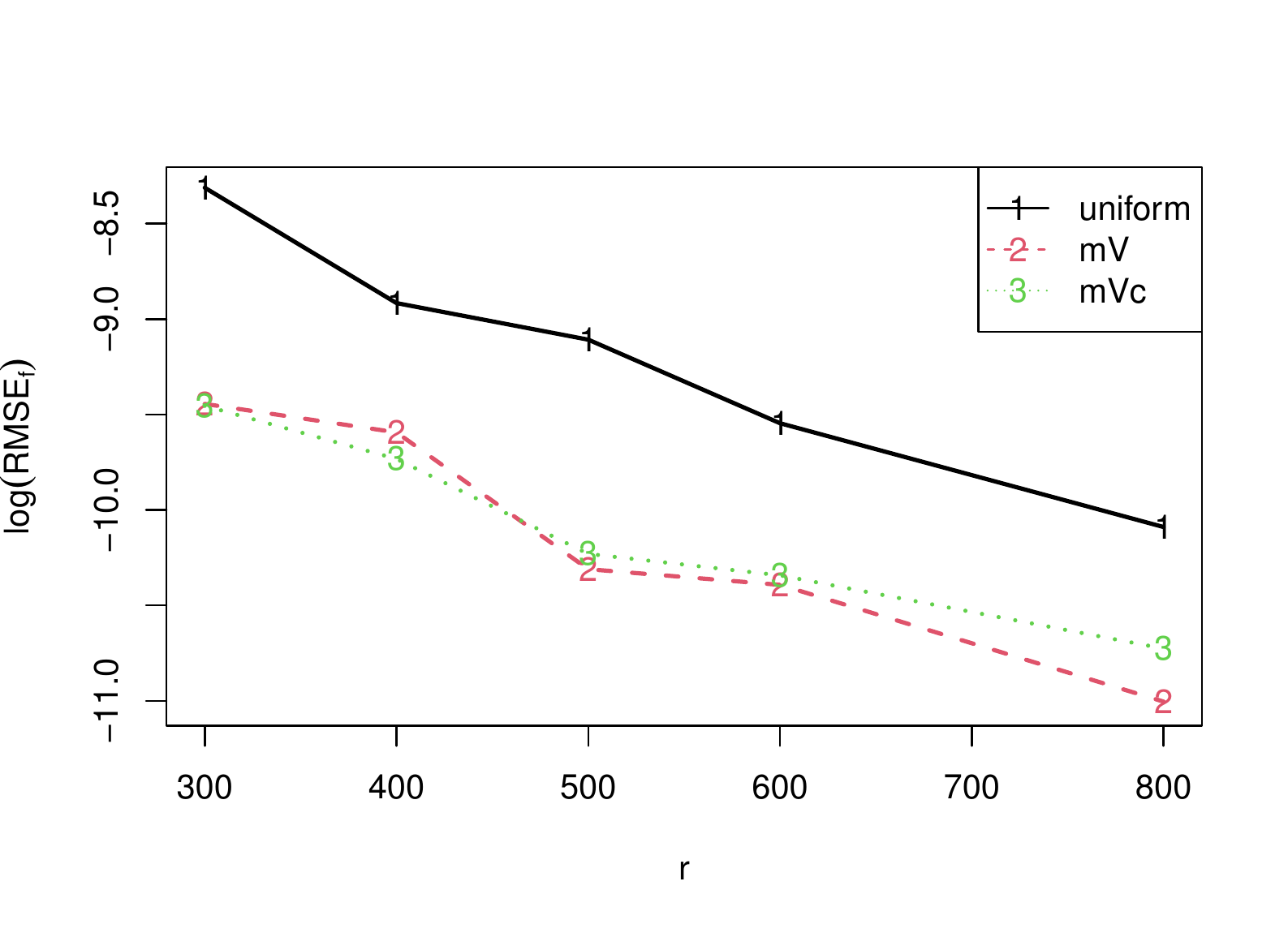}
			\caption{Accuracy of the proposed estimator  for different subsample sizes $r$ and a fixed $r_0=28$ based on mV (red), mVc (green) and uniform (black) subsampling methods.}
			\label{MSE4}
		\end{figure}
		Obviously, the estimator based on the mV and mVc criteria works better than that based on the uniform subsampling method. 
		It is worth noting that although the full sample size $n$ is not particularly large in this example, predicting the output of the highly expensive VarKarst-R model at $n$ physical observations  points is very time-consuming. The results of this example demonstrate that the proposed estimators using subsamples can approximate the full data estimates very well. It greatly saves the prediction time of the computer model as well as the computational time of the estimation for the calibration parameters.
		
		Now we construct $95\%$ confidence intervals for $a$. 
		Table \ref{Table 8} gives the results. Same as the previous example, coverage rates are not shown here also.
		\begin{table}[!ht]
			\centering
			\caption{Average lengths of $95\%$ confidence intervals for $a$.}
			\begin{tabular}{c|cccccc}
				\hline
				Criterion &   $r=300$&  $r=400$&  $r=500$&  $r=600$&  $r=800$\\
				\hline
				uniform&  $1.019$&  $0.972$&  $0.856$&  $0.685$&  $0.410$\\
				mV&  $0.586$&  $0.511$&  $0.483$&  $	0.397$&  $0.308$\\
				mVc&  $0.692$&  $0.513$&  $0.360$&  $0.335$&  $0.324$\\
				\hline
			\end{tabular}
			\label{Table 8}
		\end{table}
		We can see that the mV and mVc based subsampling methods are still better than the uniform subsampling method by comparison. 
	\end{example}

	\section{Discussion}
	\label{sec:conc}
	In this paper, we have proposed a fast calibration method with the subsamples generated by an optimal Poisson subsampling method. The consistency and asymptotic normality of the subsample-based estimators have been established. We have implemented numerical simulations and real case studies to investigate their performance. It turns out that our proposed method is effective and can extract more useful information than the uniform sampling method.


	In this work, we assume that the predictor $\hat y^s(\cdot,\cdot)$ for the computer outputs is accurate. 
	However, for non-smooth, extremely expensive computer models, it is difficult to build accurate surrogate models with limited computer experiments. In these cases, the consistency of the OLS calibration will lose. To deal with this problem, an estimator for the calibration parameter that takes into account the uncertainties of the surrogate model is needed. We leave this part for later discussion.

\bigskip
\begin{center}
{\large\bf SUPPLEMENTARY MATERIAL}
\end{center}



\begin{appendices}
	
\section{Technique proofs}

\subsection{Proof of Theorem \ref{th2:normal}}

Since $\hat{\bm\theta}-\bm\theta^*$ converges to a  normal distribution with mean zero, it is sufficient to prove that  conditional on $\mathcal D_n$ in probability,  there is 
\begin{align}
	&\tilde {\mathbf \Sigma}^{-\frac{1}{2}}(\tilde{\bm\theta}-\hat{\bm\theta})\to  N(\bm 0,\mathbf I_q) \label{eq14}
\end{align}
in distribution. Here, $\mathbf I_q$ is the identity matrix of size $q$.

Denote $\bm a=(a_1,\ldots, a_n)^T$, with $a_{i}=1$ if and only if the $i$-th data point $(\bm x_i,y_i^p)$ is included in the subsamples.	
It can be easily seen that $a_{i}$ follows from Binomial distribution $B(1,\pi_{i})$,  $i=1,2,\ldots,n$. Then there is ${\rm{E}}(a_i)=\pi_i$ and ${\rm{Var}}(a_i)=\pi_i(1-\pi_i)$.

Here we first prove the consistency of  $\tilde{\bm\theta}$.
From the definitions of  $\tilde{\bm\theta}$ and $\hat{\bm \theta}$, it suffices to prove that $l^{*}(\bm\theta)$ converges to $l(\bm\theta)$  uniformly with respect to $\bm\theta$ in conditional  probability.
According to the sampling theory, the equation (\ref{eq5}) can be represented by
\begin{equation}	
	l^{*}(\bm\theta) = \frac{1}{n}\sum_{i=1}^{n}a_{i}\frac{1}{\pi_{i}}\left[y^{p}_i-\hat y^{s}(\bm x_{i},\bm\theta)\right]^{2}.
	\notag
\end{equation}
Now, we focus on the conditional expectation and the conditional variance of $l^{*}(\bm\theta)-l(\bm\theta)$.
\begin{itemize}
	\item 	The conditional expectation of $l^{*}(\bm\theta)$  is 
	\begin{equation*}
		\label{conexp1}
		\begin{aligned}
			{\rm E}_{\bm a} [l^{*}(\bm\theta)| \mathcal{D}_{n}]
			= 
			\frac{1}{n}\sum_{i=1}^{n}\left[y^{p}_i-\hat y^{s}(\bm x_{i},\bm\theta)\right]^{2}
			=  l(\bm\theta).
		\end{aligned}
	\end{equation*} 	
	That is,
	\begin{equation}
		{\rm E}_{\bm a}\left[l^{*}(\bm\theta)-l(\bm\theta)|\mathcal{D}_{n}\right]=0.
		\label{con1}
	\end{equation}
	\item 	By some elementary calculations, the conditional variance of  $l^{*}(\bm\theta)-l(\bm\theta)$  is 
	\begin{equation}
		\begin{aligned}	\label{convar1}
			{\rm Var}_{\bm a}\left[l^{*}(\bm\theta)-l(\bm\theta)| \mathcal{D}_{n}\right]
			=& {\rm Var}_{\bm a}\left\{\frac{1}{n}\sum_{i=1}^{n}(\frac{a_{i}}{\pi_{i}}-1) \left[y^{p}_i-\hat y^{s}(\bm x_{i},\bm\theta)\right]^{2}| \mathcal{D}_{n}\right\}\\
			=& \frac{1}{n^{2}}\sum_{i=1}^{n}\frac{1-\pi_{i}}{\pi_{i}} \left[y^{p}_i-\hat y^{s}(\bm x_{i},\bm\theta)\right]^{4}\\
			\leq& O_p(\frac{1}{r}+\frac{1}{n}) \frac{1}{n}\sum_{i=1}^{n} \left[y^{p}_i-\hat y^{s}(\bm x_{i},\bm\theta)\right]^{4},
		\end{aligned}
	\end{equation}
	the last inequality is obtained from Assumption (H.\ref{a5}).
	Note that  the following equality holds from Assumption (H.\ref{a1}) and  Assumption (H.\ref{a2})
	\begin{equation}
		\sup_{\bm\theta\in\bm\Theta}	\frac{1}{n}\sum_{i=1}^{n} \left[y^{p}_i-\hat y^{s}(\bm x_{i},\bm\theta)\right]^{4} = O_p(1).
		\label{214}
	\end{equation}
	Thus, there is
	\begin{equation}
		{\rm Var}_{\bm a}\left[l^{*}(\bm\theta)-l(\bm\theta)| \mathcal{D}_{n}\right]\leq O_p(r^{-1}).
		\label{215}
	\end{equation}
\end{itemize}
The consistency of  $\tilde{\bm\theta}$ can be obtained by combining (\ref{con1}) and (\ref{215}). 

Next, we prove the asymptotic normality of  $\tilde{\bm\theta}$.
Define $\dot{l}_{j}^{*}(\bm\theta)$ as the partial derivative of $l^{*}(\bm \theta)$ with respect to $\theta_j$, $j\in\{1,\ldots,q\}$.	 
The asymptotic normality of $\tilde{\bm\theta}$ can be proved by  the Delta method for establishing asymptotic theory for M-estimation \citep{van2000asymptotic}. According to the definition of $\tilde{\bm\theta}$, we employ Taylor's theorem \citep{ferguson2017course} to expand $\dot{l}_{j}^{*}(\tilde{\bm\theta})$ at $\hat{\bm\theta}$, that is,
\begin{equation}
	\label{taylor}			0=\dot{l}_{j}^{*}(\tilde{\bm\theta})=\dot{l}_{j}^{*}(\hat{\bm\theta}) + \frac{\partial \dot{l}_{j}^{*}(\hat{\bm\theta})}{\partial \bm\theta}(\tilde{\bm\theta}-\hat{\bm\theta}) + R_j,
\end{equation}
where $R_j=(\tilde{\bm\theta}-\hat{\bm\theta})^{T}\int_{0}^{1}\int_{0}^{1}  \frac{\partial^2 \dot{l}_{j}^{*}(\hat{\bm\theta}+(\tilde{\bm\theta}-\hat{\bm\theta})uv)}{\partial \bm\theta \partial \bm\theta^T}vdudv(\tilde{\bm\theta}-\hat{\bm\theta}).$
\begin{itemize}
	\item First, we consider $\dot{l}_{j}^{*}(\hat{\bm\theta})$. Recall that $ l^{*}(\bm \theta)= \frac{1}{n}\sum_{i=1}^{r}\frac{1}{\pi_{i}}\left[y^{p}(\bm x_{i}^{*})-\hat{y}^{s}(\bm x_{i}^{*},\bm \theta)\right]^{2}$. It is obvious that, $\dot{l}_{j}^{*}(\hat{\bm\theta})$ can be represented by 
\end{itemize}
\begin{equation*}
	\dot{l}_{j}^{*}(\hat{\bm\theta})=\frac{2}{n}\sum_{i=1}^{n}\frac{a_{i}}{\pi_{i}}\left[\hat y^{s}(\bm x_{i},\hat{\bm\theta})-y^{p}_i\right]\frac{\partial \hat y^{s}(\bm x_{i},\hat{\bm\theta})}{\partial \theta_j}:=\sum_{i=1}^{n}\eta_{i}.
\end{equation*}
To find  the distribution of $\dot{l}_{j}^{*}(\hat{\bm\theta})$, the  Lindeberg-Feller condition \citep{van2000asymptotic} is verified as follows.
For every constant $\epsilon>0$,
\begin{equation}
	\begin{aligned}
		\label{LF1}
		&\sum_{i=1}^{n}{\rm E}_{\bm a}\left[ \left| \eta_{i} \right|^{2} I(\left| \eta_{i} \right|>\epsilon)|\mathcal{D}_{n} \right]\\
		\leq& \frac{1}{\epsilon}\sum_{i=1}^{n}{\rm E}_{\bm a}(|\eta_{i}|^{3}|\mathcal{D}_{n})\\
		=& \frac{1}{\epsilon}\sum_{i=1}^{n}{\rm E}_{\bm a}\left[ \frac{8a_{i}^{3}}{n^{3}\pi_{i}^{3}}\left|\hat y^{s}(\bm x_{i},\hat{\bm\theta})-y^{p}_i\right|^{3}\left| \frac{\partial \hat y^{s}(\bm x_{i},\hat{\bm\theta})}{\partial \theta_{j}} \right|^{3} | \mathcal{D}_{n} \right]\\
		\leq& \frac{1}{\epsilon}\left[ \underset{i=1,2,\ldots,n}{\max}(n\pi_{i})^{-2} \right] \frac{1}{n}\sum_{i=1}^{n}\left|\hat y^{s}(\bm x_{i},\hat{\bm\theta})-y^{p}_i\right|^{3}\left| \frac{\partial \hat y^{s}(\bm x_{i},\hat{\bm\theta})}{\partial \theta_{j}} \right|^{3} . 
	\end{aligned}
\end{equation}
By Assumptions (H.\ref{a2}) - (H.\ref{a3}), and the H\"older inequality \citep{schilling2017measures}, we have that
\begin{equation}
	\frac{1}{n}\sum_{i=1}^{n}\left|\hat{y}^{s}(\bm x_{i},\hat{\bm\theta})-y^{p}_i\right|^{3}\left| \frac{\partial y^{s}(\bm x_{i},\hat{\bm\theta})}{\partial \theta_{j}} \right|^{3}=O_p(1).
	\label{LF2}
\end{equation}

By combining (\ref{LF1}), (\ref{LF2}), and Assumption (H.\ref{a5}), we obtain
\begin{equation*}
		\sum_{i=1}^{n}{\rm E}_{\bm a}\left[ \left| \eta_{i} \right|^{2} I(\left| \eta_{i} \right|>\epsilon)|\mathcal{D}_{n} \right] 
		\leq \frac{1}{\epsilon}O_p(\frac{1}{r^{2}})O_p(1)=o_p(1).
		\label{LF} 
\end{equation*}
Thus, by the Lindeberg-Feller central limit theorem \citep{van2000asymptotic}, it can be concluded that as $n\to\infty$ and $r\to \infty$, conditional on $\mathcal D_n$,
\begin{equation}\label{dl}
	\dot{l}_{j}^{*}(\hat{\bm\theta})|\mathcal D_n \sim  N\left({\rm E}_{\bm a} [\dot{l}_{j}^{*}(\hat{\bm\theta})|\mathcal{D}_{n}],{\rm Var}_{\bm a}[\dot{l}_{j}^{*}(\hat{\bm\theta})|\mathcal{D}_{n}]\right).
\end{equation}	

By some simple calculations, the conditional expectation of $\dot{l}_{j}^{*}(\hat{\bm\theta})$ is
\begin{equation*}
	\begin{aligned}
		{\rm E}_{\bm a}\left[ \dot{l}_{j}^{*}(\hat{\bm\theta}) |\mathcal{D}_{n}\right]
		=\frac{2}{n}\sum_{i=1}^{n} \left[\hat y^{s}(\bm x_{i},\hat{\bm\theta})-y^p_i\right] \frac{\partial \hat y^{s}(\bm x_{i},\hat{\bm\theta})}{\partial \theta_j},   \\
	\end{aligned}
\end{equation*}	
and  the conditional variance of  $\dot{l}_{j}^{*}(\hat{\bm\theta})$ can be evaluated by
\begin{equation}
	\begin{aligned}
		{\rm Var}_{\bm a}\left[ \dot{l}_{j}^{*}(\hat{\bm\theta})|\mathcal{D}_{n} \right] 
		=\frac{4}{n^{2}}\sum_{i=1}^{n}\frac{1-\pi_{i}}{\pi_{i}}\left\{\left[\hat y^{s}(\bm x_{i},\hat{\bm\theta})-y^{p}(\bm x_{i})\right]\frac{\partial \hat y^{s}(\bm x_{i},\hat{\bm\theta})}{\partial\theta_j}\right\}^{2}:=\tilde{\mathbf V}_{j}, 
		\notag
	\end{aligned}
\end{equation}
so it shows that $\dot{l}_{j}^{*}(\hat{\bm\theta})=O_{p|\mathcal{D}_{n}}(r^{-1/2})$.

According to the define of $\hat{\bm\theta}$ 	(\ref{eq3-expensive}), we have that ${\rm E}_{\bm a}[\dot{l}_{j}^{*}(\hat{\bm\theta})|\mathcal{D}_{n}]=0$.
By plugging the conditional expectation and the conditional variance of  $\dot{l}_{j}^{*}(\hat{\bm\theta})$ into  (\ref{dl}), we have that as $n\to\infty$ and $r\to \infty$, conditional on $\mathcal D_n$,
\begin{equation}\label{dl-2}
	\dot{l}_j^{*}(\hat{\bm\theta})\rightarrow N(0,\tilde{\mathbf V}_{j})
\end{equation}	
in distribution.	

\begin{itemize}
	\item Next, we consider $\frac{\partial \dot{l}_{j}^{*}(\hat{\bm\theta})}{\partial \bm\theta^T}$, which can be written as
\end{itemize}
\begin{equation*} \label{ddl}
	\begin{aligned}	
		\frac{\partial \dot{l}_j^{*}(\hat{\bm\theta})}{\partial \bm\theta^T}
		&=\frac{1}{n}\frac{\partial^{2}\sum_{i=1}^{n}\frac{a_{i}}{\pi_{i}}\left[y^p_i-\hat y^{s}(\bm x_{i},\hat{\bm\theta})\right]^{2}}{\partial\bm\theta^T \partial \theta_j}.
	\end{aligned}
\end{equation*}		
By the Lebesgue's dominated convergence theorem,  the conditional expectation of $\frac{\partial \dot{l}_{j}^{*}(\hat{\bm\theta})}{\partial \bm\theta^T}$ can be represented by
\begin{equation*}
	\begin{aligned}	
		{\rm{E}}_{\bm a }	\left[\frac{\partial \dot{l}_j^{*}(\hat{\bm\theta})}{\partial \bm\theta^T}\right]
		= \frac{1}{n}\sum_{i=1}^{n} \frac{\partial^{2}\left[y^p_i-\hat y^{s}(\bm x_i,\hat{\bm\theta})\right]^{2}}{\partial\bm\theta^T \partial\theta_j},
		\label{ddlJ}
	\end{aligned}			
\end{equation*}
and the conditional variance of  $\frac{\partial \dot{l}_{j}^{*}(\hat{\bm\theta})}{\partial \bm\theta^T}$ is
\begin{equation*}
	{\rm{Var}}_{\bm a }	\left[\frac{\partial \dot{l}_j^{*}(\hat{\bm\theta})}{\partial \bm\theta^T}\right]
	= \frac{1}{n^2}\sum_{i=1}^{n} \frac{1-\pi_i}{\pi_i} \frac{\partial^{2}
		\left[y^p_i-\hat y^{s}(\bm x_i,\hat{\bm\theta})\right]^{2}}{\partial\bm\theta^T \partial\theta_j}	
\end{equation*}
By using the similar arguments in proving (\ref{215}), we have that under Assumptions (H.\ref{a2}) – (H.\ref{a5}), ${\rm{Var}}_{\bm a }	\left[\frac{\partial \dot{l}_j^{*}(\hat{\bm\theta})}{\partial \bm\theta^T}\right]=O_p(r^{-1})$. As a result,
\begin{equation}
	\label{515}
	\frac{\partial \dot{l}_j^{*}(\hat{\bm\theta})}{\partial \bm\theta^T}=\frac{1}{n}\sum_{i=1}^{n} \frac{\partial^{2}\left[y^p_i-\hat y^{s}(\bm x_i,\hat{\bm\theta})\right]^{2}}{\partial\bm\theta^T \partial\theta_j}+O_{p|\mathcal{D}_{n}}(r^{-1/2})= \tilde{\mathbf{J}}+O_{p|\mathcal{D}_{n}}(r^{-1/2}).
\end{equation}

\begin{itemize}	
	\item  Third, we bound $\int_{0}^{1}\int_{0}^{1}  \frac{\partial^2 \dot{l}_{j}^{*}(\hat{\bm\theta}+(\tilde{\bm\theta}-\hat{\bm\theta})uv)}{\partial \bm\theta \partial \bm\theta^T}vdudv$. 
\end{itemize}
Suppose $\bm\theta_0\in \bm\Theta_0$, next we bound 
$ \frac{\partial^2 \dot{l}_{j}^{*}({\bm\theta_0})}{\partial \theta_k \partial \theta_{k'}}$, which can be written as follows
\begin{equation*} 
	\begin{aligned}	
		\frac{\partial \dot{l}_j^{*}({\bm\theta}_0)}{\partial \theta_k\partial \theta_{k'}}
		&=\frac{1}{n}\frac{\partial^{3}\sum_{i=1}^{n}\frac{a_{i}}{\pi_{i}}\left[y^p_i-\hat y^{s}(\bm x_{i},{\bm\theta}_0)\right]^{2}}{\partial\theta_k\partial \theta_{k'} \partial \theta_j}.
	\end{aligned}
\end{equation*}		
Similar with (\ref{515}), $\frac{\partial \dot{l}_j^{*}({\bm\theta}_0)}{\partial \theta_k\partial \theta_{k'}}$ can be bounded by
$$\frac{1}{n}\frac{\partial^{3}\sum_{i=1}^{n}\left[y^p_i-\hat y^{s}(\bm x_{i},{\bm\theta}_0)\right]^{2}}{\partial\theta_k\partial \theta_{k'} \partial \theta_j}+O_p(r^{-1/2}).$$

Because $\|\tilde{\bm\theta}-\hat{\bm\theta}\|=o_{p|\mathcal{D}_{n}}(1)$, $\hat{\bm\theta}+(\tilde{\bm\theta}-\hat{\bm\theta})uv\in \bm\Theta_0$, it follows that
\begin{equation*}
	\label{dddl}
	\int_{0}^{1}\int_{0}^{1}  \frac{\partial^2 \dot{l}_{j}^{*}(\hat{\bm\theta}+(\tilde{\bm\theta}-\hat{\bm\theta})uv)}{\partial \theta_k \partial \theta_{k'}}vdudv=\int_{0}^{1}\int_{0}^{1} \frac{1}{n}
	\frac{\partial^{3}\sum_{i=1}^{n}\left[y^p_i-\hat y^{s}(\bm x_{i},{\bm\theta}_0)\right]^{2}}{\partial\theta_k\partial \theta_{k'} \partial \theta_j}vdudv+O_p(r^{-\frac{1}{2}})
\end{equation*}
Because $\frac{1}{n}
\frac{\partial^{3}\sum_{i=1}^{n}\left[y^p_i-\hat y^{s}(\bm x_{i},{\bm\theta}_0)\right]^{2}}{\partial\theta_k\partial \theta_{k'} \partial \theta_j}=O_p(1)$, we have
\begin{equation}
	\label{dddl3}
	\int_{0}^{1}\int_{0}^{1}  \frac{\partial^2 \dot{l}_{j}^{*}(\hat{\bm\theta}+(\tilde{\bm\theta}-\hat{\bm\theta})uv)}{\partial \theta_k \partial \theta_{k'}}vdudv=O_p(1).
\end{equation}

By plugging (\ref{dl-2}), (\ref{515}), (\ref{dddl3}) into (\ref{taylor}), there is
\begin{equation*}
	\begin{aligned}
	\tilde{\mathbf{\Sigma}}^{-1/2}(\tilde{\bm\theta}-\hat{\bm\theta}) &= -\tilde{\mathbf{\Sigma}}^{-1/2}\left[\frac{\partial \dot{l}^{*}(\hat{\bm\theta})}{\partial \bm\theta^T}\right]^{-1}\dot{l}^{*}(\hat{\bm\theta}) + O_p(r^{-1/2})\\
		&= -\tilde{\mathbf{\Sigma}}^{-1/2}\tilde{\mathbf{J}}^{-1}\dot{l}^{*}(\hat{\bm\theta}) - \tilde{\mathbf{\Sigma}}^{-1/2}\left\{\left[\frac{\partial \dot{l}^{*}(\hat{\bm\theta})}{\partial \bm\theta^T}\right]^{-1} - \tilde{\mathbf{J}}^{-1} \right\}\dot{l}^{*}(\hat{\bm\theta}) + O_p(r^{-1/2})\\
		&=-\tilde{\mathbf{\Sigma}}^{-1/2}\tilde{\mathbf{J}}^{-1}\tilde{\mathbf{V}}^{1/2}\tilde{\mathbf{V}}^{-1/2}\dot{l}^{*}(\hat{\bm\theta}) + O_p(r^{-1/2}).
	\end{aligned}
\end{equation*}
The desired result  follows from the Slutsky's Theorem \citep{ferguson2017course} and the fact that 
\begin{align*}
\tilde{\mathbf{\Sigma}}^{-1/2}\tilde{\mathbf{J}}^{-1}\tilde{\mathbf{V}}^{1/2}(\tilde{\mathbf{\Sigma}}^{-1/2}\tilde{\mathbf{J}}^{-1}\tilde{\mathbf{V}}^{1/2})^T=\mathbf{I_q}.
\end{align*}


\subsection{Proof of Theorem \ref{th3:pi}}

The subsampling probabilities can be obtained by minimizing the following objective function:
\begin{equation}
	\begin{aligned} \label{tr(V)}
		\min& \quad {\rm tr}\left\{ \tilde{\mathbf J}^{-1} \frac{4}{n^{2}}\sum_{i=1}^{n}\frac{1}{\pi_{i}} \left[y^p_i-\hat y^{s}(\bm x_{i},\hat{\bm\theta})\right]^{2} \frac{\partial \hat y^{s}(\bm x_{i},\hat{\bm\theta})}{\partial\bm\theta^{T}}\frac{\partial\hat y^{s}(\bm x_{i},\hat{\bm\theta})}{\partial\bm\theta} \tilde{\mathbf J}^{-1}  \right\}\\
		s.t.& \quad \sum_{i=1}^{n}\pi_{i}=r, \, 0\leq\pi_{i}\leq1, \, i=1,\ldots,n.
	\end{aligned}	
\end{equation}
Based on the property of the trace of a product, we have that
\begin{equation}
	\label{eq3.3}
	\begin{aligned}
		&{\rm tr}\left\{ \tilde{\mathbf J}^{-1} \frac{4}{n^{2}}\sum_{i=1}^{n}\frac{1}{\pi_{i}} \left[y^p_i-\hat y^{s}(\bm x_{i},\hat{\bm\theta})\right]^{2} \frac{\partial \hat y^{s}(\bm x_{i},\hat{\bm\theta})}{\partial\bm\theta^{T}}\frac{\partial \hat y^{s}(\bm x_{i},\hat{\bm\theta})}{\partial\bm\theta} \tilde{\mathbf J}^{-1}  \right\}\\
		=&\frac{4}{n^{2}}\sum_{i=1}^{n}\frac{1}{\pi_{i}}  \left[y^p_i-\hat y^{s}(\bm x_{i},\hat{\bm\theta})\right]^{2} \frac{\partial \hat y^{s}(\bm x_{i},\hat{\bm\theta})}{\partial\bm\theta} \tilde{\mathbf J}^{-2}\frac{\partial \hat y^{s}(\bm x_{i},\hat{\bm\theta})}{\partial\bm\theta^{T}}\\
		=&\frac{4}{n^{2}}\left(\frac{1}{r}\sum_{i=1}^{n}\pi_{i}\right)\sum_{i=1}^{n}\frac{1}{\pi_{i}} \left[y^p_i-\hat y^{s}(\bm x_{i},\hat{\bm\theta})\right]^{2} \frac{\partial \hat y^{s}(\bm x_{i},\hat{\bm\theta})}{\partial\bm\theta} \tilde{\mathbf J}^{-2}\frac{\partial\hat  y^{s}(\bm x_{i},\hat{\bm\theta})}{\partial\bm\theta^{T}}  ,
	\end{aligned}
\end{equation}
the last equality is due to the fact that $\sum_{i=1}^{n}\pi_{i}=r$.

Define $h_{i}^{mV}=\left| y^p_i-\hat y^{s}(\bm x_{i},\hat{\bm\theta})\right|   \left[ \frac{\partial \hat y^{s}(\bm x_{i},\hat{\bm\theta})}{\partial\bm\theta} \tilde{\mathbf J}^{-2}\frac{\partial \hat y^{s}(\bm x_{i},\hat{\bm\theta})}{\partial\bm\theta^{T}} \right]^{1/2}, i=1,\ldots,n$.
Without losing generality, we assume $h_{1}^{mV}\leq h_{2}^{mV}\leq \ldots \leq h_{n}^{mV} $.  (\ref{eq3.3}) can continue to be represented by 
\begin{equation}
	\begin{aligned}
		&{\rm tr}\left\{ \tilde{\mathbf J}^{-1} \frac{4}{n^{2}}\sum_{i=1}^{n}\frac{1}{\pi_{i}} \left[y^p_i-\hat y^{s}(\bm x_{i},\hat{\bm\theta})\right]^{2} \frac{\partial \hat y^{s}(\bm x_{i},\hat{\bm\theta})}{\partial\bm\theta^{T}}\frac{\partial \hat y^{s}(\bm x_{i},\hat{\bm\theta})}{\partial\bm\theta} \tilde{\mathbf J}^{-1}  \right\}\\
		\geq& \frac{4}{n^{2}r}\left\{ \sum_{i=1}^{n} \left|y^p_i-\hat y^{s}(\bm x_{i},\hat{\bm\theta})\right|  \left[\frac{\partial \hat y^{s}(\bm x_{i},\hat{\bm\theta})}{\partial\bm\theta^{T}} \tilde{\mathbf J}^{-2}\frac{\partial \hat y^{s}(\bm x_{i},\hat{\bm\theta})}{\partial\bm\theta}   \right] ^{\frac{1}{2}}  \right\}^{2}\\
		=&\frac{4}{n^{2}r} \left( \sum_{i=1}^{n} h_{i}^{mV} \right)^{2},
		\notag
	\end{aligned}
\end{equation}
the inequality is from the Cauchy-Schwarz inequality and the equality holds if and only if $\pi_{i} \propto h_{i}^{mV}$. Next, we consider the following two cases:
\begin{itemize}
	\item  For $i=1,\ldots,n$,  $\pi_{i}^{mV} = r\frac{h_{i}^{mV}}{\sum_{j=1}^{n} h_{j}^{mV}} \leq 1$, then $\left\{ \pi_{i}^{mV} \right\}_{i=1}^{n}$ give the optimal solution.
	\item If there exists some $i$ such that $r\frac{h_{i}^{mV} }{\sum_{j=1}^{n} h_{j}^{mV}} > 1$, and from the definition of $k$, we know that the number of such $i$ is $k$, then the initial minimizing formula (\ref{tr(V)}) is transformed into the following problem:
\end{itemize}
\begin{equation}
	\begin{aligned} \label{ttttr(V)}
		\min& \quad {\rm tr}\left\{ \tilde{\mathbf J}^{-1} \frac{4}{n^{2}}\sum_{i=1}^{n-k}\frac{1}{\pi_{i}} \left[y^p_i-y^{s}(\bm x_{i},\hat{\bm\theta})\right]^{2} \frac{\partial y^{s}(\bm x_{i},\hat{\bm\theta})}{\partial\bm\theta}\frac{\partial y^{s}(\bm x_{i},\hat{\bm\theta})}{\partial\bm\theta^{T}} \tilde{\mathbf J}^{-1}  \right\}\\
		s.t.& \quad \sum_{i=1}^{n-k}\pi_{i}=r-k, \, 0\leq\pi_{i}\leq1, \, i=1,\ldots,n-k,\\
		& \quad \pi_{n-k+1},\ldots,\pi_{n}=1.
	\end{aligned}	
\end{equation}

Similar to (\ref{eq3.3}), by applying the Cauchy-Schwarz inequality, it holds that
\begin{equation}
\label{besttr}
	\begin{aligned}	
		&{\rm tr}\left\{ \tilde{\mathbf J}^{-1} \frac{4}{n^{2}}\sum_{i=1}^{n-k}\frac{1}{\pi_{i}} \left[y^p_i-y^{s}(\bm x_{i},\hat{\bm\theta})\right]^{2} \frac{\partial y^{s}(\bm x_{i},\hat{\bm\theta})}{\partial\bm\theta}\frac{\partial y^{s}(\bm x_{i},\hat{\bm\theta})}{\partial\bm\theta^{T}} \tilde{\mathbf J}^{-1}  \right\}\\		
		\geq& \frac{4}{n^{2}(r-k)}\left\{  \sum_{i=1}^{n-k}  \left|y^p_i-y^{s}(\bm x_{i},\hat{\bm\theta})\right|   \left[\frac{\partial y^{s}(\bm x_{i},\hat{\bm\theta})}{\partial\bm\theta^{T}} \tilde{\mathbf J}^{-2}\frac{\partial y^{s}(\bm x_{i},\hat{\bm\theta})}{\partial\bm\theta}   \right] ^{\frac{1}{2}}  \right\}^{2}\\
		=&\frac{4}{n^{2}(r-k)}\left( 	\sum_{i=1}^{n-k}h_{i}^{mV} \right) ^{2},	
	\end{aligned}
\end{equation}
the equality in the last inequality holds if and only if $\pi_{i} \propto h_{i}^{mV}$.

Suppose there exists $M$ such that
\begin{equation}\label{defM}
\underset{i=1,\ldots,n}{\max}\frac{h_{i}^{mV}\wedge M}{\sum_{j=1}^{n} (h_{j}^{mV}\wedge M)}=\frac{1}{r}\end{equation}
and $h_{n-k}^{mV}<M\leq h_{n-k+1}^{mV}$. Let $\pi_{i}^{mV}=r\frac{h_{i}^{mV}\wedge M}{\sum_{j=1}^{n} (h_{j}^{mV}\wedge M)}$. Next, we verify that $\pi_{i}^{mV}$  is the optimal probability. 
Combining  the fact that $\max_{i=1,\ldots,n} h_i^{mV}\wedge M =M$ and (\ref{defM}), there is  $\sum_{i=1}^{n-k}h_{i}^{mV}=(r-k)M$. Thus, $\pi_{i}^{mV}=(r-k)\frac{h_{i}^{mV}}{\sum_{j=1}^{n-k} h_{j}^{mV}}=\frac{h_{i}^{mV}}{M},i=1,\ldots,n-k$ and $\pi_{i}^{mV}=1, i=n-k+1,\ldots,n$.  From the  Cauchy-Schwarz inequality, it is clear that $\pi_{i}^{mV}$ is the optimal probability that makes the last equality hold in (\ref{besttr}).

Now we will prove the existence of $M$ and prove that the range of $M$ is $h_{n-k}^{mV}<M\leq h_{n-k+1}^{mV}$. From the definition of $k$, we have that	
\begin{equation*}
	\frac{(r-k+1)h_{n-k+1}^{mV}}{\sum_{i=1}^{n-k+1}h_{i}^{mV}} \geq 1 
	\quad \mbox{and} \quad
	\frac{(r-k)h_{n-k}^{mV}}{\sum_{i=1}^{n-k}h_{i}^{mV}} < 1.
\end{equation*}

Denote $M_{1}=h_{n-k+1}^{mV}$, $M_{2}=h_{n-k}^{mV}$, then there is 
\begin{equation*}
	\frac{(r-k+1)h_{n-k+1}^{mV}+(k-1)M_{1}}{\sum_{i=1}^{n-k+1}h_{i}^{mV}+(k-1)M_{1}} \geq 1 
	\quad \mbox{and} \quad
	\frac{(r-k)h_{n-k}^{mV}+kM_{2}}{\sum_{i=1}^{n-k}h_{i}^{mV}+kM_{2}} < 1,
\end{equation*}
that is $r\frac{h_{i}^{mV}\wedge M_{1}}{\sum_{j=1}^{n} (h_{j}^{mV}\wedge M_{1})} \geq 1$ and  $r\frac{h_{i}^{mV}\wedge M_{2}}{\sum_{j=1}^{n} (h_{j}^{mV}\wedge M_{2})} < 1$. 
Thus the existence of $M$ and $M_{2}<M\leq M_{1}$ are attributed to the continuity of $\underset{i=1,\ldots,n}{\max}\frac{h_{i}^{mV}\wedge M}{\sum_{j=1}^{n} (h_{j}^{mV}\wedge M)}$. 

Otherwise, $\forall$ $ h_{n}^{mV}\geq M^{'}>M $, $M^{'}\wedge h_{n}^{mV} \geq M\wedge h_{n}^{mV}$, and $\frac{M^{'}}{M}\sum_{i=1}^{n} h_{i}^{mV}\wedge M \geq \sum_{i=1}^{n} h_{i}^{mV}\wedge M^{'}$, so $\frac{h_{n}^{mV}\wedge M}{\sum_{j=1}^{n} (h_{j}^{mV}\wedge M)}$ is nondecreasing on $M\in (h_{1}^{mV},h_{n}^{mV})$. Therefore
$$\underset{i=1,\ldots,n}{\max}\frac{h_{i}^{mV}\wedge M}{\sum_{j=1}^{n} (h_{j}^{mV}\wedge M)}=\frac{1}{r}.$$
It indicates that $h_{n-k}^{mV}<M\leq h_{n-k+1}^{mV}$, the proof is completed.

\subsection{Proof of Theorem \ref{th4:pi}}

The proof is similar to the proof to Theorem \ref{th3:pi}, so we ignore details.

\subsection{Proof of Theorem \ref{th4:cvg}}

Since $\breve{\pi}_{i}^{w} \geq \rho r/n$, we have $\underset{i=1,\ldots,n}{\max}(n\pi_i)^{-1}=O_p(r^{-1})$, Theorem \ref{th2:normal} indicates Theorem \ref{th4:cvg}.

\subsection{Proof of Theorem \ref{th6:normal}}

Similar with the proof of Theorem \ref{th2:normal}, it is sufficient to prove that  conditional on $\mathcal D_n$ in probability,  
\begin{align}
		&\breve{\mathbf \Sigma }^{-\frac{1}{2}}(\breve{\bm\theta}-\hat{\bm\theta}) \to  N(\bm 0,\mathbf I_q)		
	\end{align} 
	in distribution.

Recall that $ l_{\tilde{\bm\theta}_{0}}^{*}(\bm\theta)=\frac{1}{n}\sum_{i=1}^{r}\frac{1}{\breve{\pi}_{i}^{w}\wedge 1}\left[y^{p}(\bm x_{i}^{*})-\hat y^{s}(\bm x_{i}^{*},\bm \theta)\right]^{2} $. 
The asymptotic normality of $\breve{\bm\theta}$ can be proved by  the Delta method for establishing asymptotic theory for M-estimation \citep{van2000asymptotic}. According to the definition of $\breve{\bm\theta}$, we use Taylor's theorem \citep{ferguson2017course} to expand $\frac{\partial l_{\tilde{\bm\theta}_{0}}^{*}(\breve{\bm\theta})}{\partial \theta_j}$ at $\hat{\bm\theta}$, that is,
\begin{equation}
	\label{taylor2}
	0=\frac{\partial l_{\tilde{\bm\theta}_{0}}^{*}(\breve{\bm\theta})}{\partial \theta_j}=\frac{\partial l_{\tilde{\bm\theta}_{0}}^{*}(\hat{\bm\theta})}{\partial \theta_j} + \frac{\partial^2 l_{\tilde{\bm\theta}_{0}}^{*}(\hat{\bm\theta})}{\partial\bm\theta \partial \theta_j}(\breve{\bm\theta}-\hat{\bm\theta}) + R_j,
\end{equation}
where $R_j=(\breve{\bm\theta}-\hat{\bm\theta})^{T}\int_{0}^{1}\int_{0}^{1}  
\frac{\partial^3 l_{\tilde{\bm\theta}_{0}}^{*}(\hat{\bm\theta}+(\breve{\bm\theta}-\hat{\bm\theta})uv)}{\partial\bm\theta \partial\bm\theta^T \partial\theta_j}vdudv(\breve{\bm\theta}-\hat{\bm\theta})$. 

\begin{itemize}
	\item First, we consider $\frac{\partial l_{\tilde{\bm\theta}_{0}}^{*}(\hat{\bm\theta})}{\partial \theta_j}$, which can be represented by
\end{itemize}
\begin{equation*}
	\frac{\partial l_{\tilde{\bm\theta}_{0}}^{*}(\hat{\bm\theta})}{\partial \theta_j}=\frac{2}{n}\sum_{i=1}^{n}a_{i}\frac{1}{\breve{\pi}_i^w\wedge 1}\left[\hat y^{s}(\bm x_{i},\hat{\bm\theta})-y^{p}_i\right]\frac{\partial \hat y^{s}(\bm x_{i},\hat{\bm\theta})}{\partial \theta_j}:=\sum_{i=1}^{n}\tilde{\eta}_{i}.
\end{equation*}
To find distribution of $\frac{\partial l_{\tilde{\bm\theta}_{0}}^{*}(\hat{\bm\theta})}{\partial \theta_j}$, the Lindeberg-Feller condition \citep{van2000asymptotic} is verified as follows.
For every constant $\epsilon>0$,
\begin{equation*}
	\begin{aligned}
		&\sum_{i=1}^{n}{\rm E}_{\bm a}\left[ \left| \tilde{\eta}_{i} \right|^{2} I(\left| \tilde{\eta}_{i} \right|>\epsilon)|\mathcal{D}_{n} \right]\\
		\leq& \frac{1}{\epsilon}\sum_{i=1}^{n}{\rm E}_{\bm a}\left[ \frac{8a_{i}^{3}}{n^{3}(\breve{\pi}_i^w\wedge 1)^3}\left|\hat y^{s}(\bm x_{i},\hat{\bm\theta})-y^{p}_i\right|^{3}\left| \frac{\partial \hat y^{s}(\bm x_{i},\hat{\bm\theta})}{\partial \theta_{j}} \right|^{3} |\mathcal{D}_{n} \right]\\
		=&o_{p}(1),
	\end{aligned}
\end{equation*}
the last equality holds for a similar reason to (\ref{LF}).

Thus, by the Lindeberg-Feller central limit theorem \citep{van2000asymptotic}, it can be concluded that as $n\to\infty$ and $r\to \infty$, conditional on $\mathcal D_n$,
\begin{equation}\label{dl_2}
	\frac{\partial l_{\tilde{\bm\theta}_{0}}^{*}(\hat{\bm\theta})}{\partial \theta_j} \rightarrow 	N\left({\rm E}_{\bm a} \left[\frac{\partial l_{\tilde{\bm\theta}_{0}}^{*}(\hat{\bm\theta})}{\partial \theta_j} |\mathcal{D}_{n}\right],{\rm Var}_{\bm a}\left[\frac{\partial l_{\tilde{\bm\theta}_{0}}^{*}(\hat{\bm\theta})}{\partial \theta_j}|\mathcal{D}_{n}\right]\right)
\end{equation}	
in distribution.

By some simple calculations, the conditional expectation of $\frac{\partial l_{\tilde{\bm\theta}_{0}}^{*}(\hat{\bm\theta})}{\partial \theta_j}$ is
\begin{equation*}
	\begin{aligned}
		{\rm E}_{\bm a}\left[ \frac{\partial l_{\tilde{\bm\theta}_{0}}^{*}(\hat{\bm\theta})}{\partial \theta_j} |\mathcal{D}_{n}\right]
		=\frac{2}{n}\sum_{i=1}^{n} \left[\hat{y}^{s}(\bm x_{i},\hat{\bm\theta})-y^p_i\right] \frac{\partial \hat{y}^{s}(\bm x_{i},\hat{\bm\theta})}{\partial \theta_j}, 
	\end{aligned}
\end{equation*}	
and the conditional variance of $\frac{\partial l_{\tilde{\bm\theta}_{0}}^{*}(\hat{\bm\theta})}{\partial \theta_j}$ can be evaluated by
\begin{equation}
	{\rm Var}_{\bm a}\left[ \frac{\partial l_{\tilde{\bm\theta}_{0}}^{*}(\hat{\bm\theta})}{\partial \theta_j}|\mathcal{D}_{n} \right] 
	=\frac{4}{n^{2}}\sum_{i=1}^{n}\frac{1-\breve{\pi}_i^w\wedge 1}{\breve{\pi}_i^w\wedge 1}\left\{\left[\hat{y}^{s}(\bm x_{i},\hat{\bm\theta})-y_i^{p}\right]\frac{\partial \hat{y}^{s}(\bm x_{i},\hat{\bm\theta})}{\partial\theta_j}\right\}^{2}:=\breve	{\mathbf V}_j^w.
	\notag
\end{equation}

According to the define of $\hat{\bm\theta}$, we have that ${\rm E}_{\bm a}\left[ \frac{\partial l_{\tilde{\bm\theta}_{0}}^{*}(\hat{\bm\theta})}{\partial \theta_j} |\mathcal{D}_{n}\right]=0$.
By plugging the conditional expectation and the conditional variance of $\frac{\partial l_{\tilde{\bm\theta}_{0}}^{*}(\hat{\bm\theta})}{\partial \theta_j}$ into (\ref{dl_2}), we have that as $n\to\infty$ and $r\to \infty$, conditional on $\mathcal D_n$,
\begin{equation}\label{dl2}
	\frac{\partial l_{\tilde{\bm\theta}_{0}}^{*}(\hat{\bm\theta})}{\partial \theta_j} {\rightarrow} N(0,\breve	{\mathbf V}_j^w).
\end{equation}	
in distribution.

Now, we discuss the distance between $\breve	{\mathbf V}^w$ and $\tilde	{\mathbf V}^w$. Let $\|\bm A\|_s$ be the spectral norm of a vector or matrix $\bm A$, there is
\begin{equation*}
	\begin{aligned}
		&\left\|  \breve	{\mathbf V}^w-\tilde	{\mathbf V}^w \right\|_s \\
		=&\left\|  \frac{4}{n^{2}}\sum_{i=1}^{n}
		\frac{1}{\pi_i^w\wedge 1} \left(\frac{\pi_i^w\wedge 1}{\breve{\pi}_i^w\wedge 1}-1\right)
		\left[\hat y^{s}(\bm x_{i},\hat{\bm\theta})-y_i^{p}\right]^{2}\frac{\partial \hat y^{s}(\bm x_{i},\hat{\bm\theta})}{\partial\bm\theta^T}\frac{\partial \hat y^{s}(\bm x_{i},\hat{\bm\theta})}{\partial\bm\theta} \right\|_s \\
			\leq&\frac{4}{n^{2}}\sum_{i=1}^{n}
		\frac{1}{\pi_i^w\wedge 1}
		\left| \frac{\pi_i^w\wedge 1}{\breve{\pi}_i^w\wedge 1}-1\right| \left[\hat y^{s}(\bm x_{i},\hat{\bm\theta})-y_i^{p}\right]^{2}   \frac{\partial \hat y^{s}(\bm x_{i},\hat{\bm\theta})}{\partial\bm\theta} \frac{\partial \hat y^{s}(\bm x_{i},\hat{\bm\theta})}{\partial\bm\theta^T} \\
		\leq&\left( \operatorname*{max}_{i=1,\ldots,n}\frac{1}{n\pi_i^w}\right)\frac{4}{n}\sum_{i=1}^{n}\left| \frac{\pi_i^w\wedge 1}{\breve{\pi}_i^w\wedge 1}-1\right| \left[\hat y^{s}(\bm x_{i},\hat{\bm\theta})-y_i^{p}\right]^{2}   \frac{\partial \hat y^{s}(\bm x_{i},\hat{\bm\theta})}{\partial\bm\theta} \frac{\partial \hat y^{s}(\bm x_{i},\hat{\bm\theta})}{\partial\bm\theta^T}   \\
		\leq&\frac{1}{\rho r}\frac{4}{n}\sum_{i=1}^{n}\left| \frac{\pi_i^w\wedge 1}{\breve{\pi}_i^w\wedge 1}-1\right|\left[\hat y^{s}(\bm x_{i},\hat{\bm\theta})-y_i^{p}\right]^{2} \frac{\partial \hat y^{s}(\bm x_{i},\hat{\bm\theta})}{\partial\bm\theta} \frac{\partial \hat y^{s}(\bm x_{i},\hat{\bm\theta})}{\partial\bm\theta^T},
	\end{aligned}
\end{equation*}
where the last inequality holds from (\ref{pi weight}). We consider the case under the mV criterion, and the same goes for the mVc criterion. 
Let $h_0^{mV}(\bm x_i)$ be the value of ${h}_i^{mV}$ with $\hat{\bm\theta}$  replaced by the pilot estimator $\tilde{\bm\theta}_0$. Since 
\begin{equation*}
	\begin{aligned}
		\left| \frac{\pi_i^w\wedge 1}{\breve{\pi}_i^w\wedge 1}-1\right| 
		\leq& \left| \frac{\pi_i^w-\breve{\pi}_i^w}{\breve{\pi}_i^w} \right| 
		\leq (1-\rho)\left| \frac{\frac{h_0^{mV}(\bm x_i)}{\sum_{i=1}^{n}h_0^{mV}(\bm x_i)}-\frac{h_i^{mV}}{\sum_{i=1}^{n}h_i^{mV}}}{\rho r\frac{1}{n}} \right|   \\
		\leq&  \frac{ \left| \frac{h_0^{mV}(\bm x_i)}{\sum_{i=1}^{n}h_0^{mV}(\bm x_i)}-\frac{h_i^{mV}}{\sum_{i=1}^{n}h_0^{mV}(\bm x_i)} \right|  + \left| \frac{h_i^{mV}}{\sum_{i=1}^{n}h_0^{mV}(\bm x_i)}- \frac{h_i^{mV}}{\sum_{i=1}^{n}h_i^{mV}}\right| }{\rho r\frac{1}{n}},
	\end{aligned}
\end{equation*}
by some simple calculations, we have that 
\begin{equation}
	\label{disV1}
	\begin{aligned}
		&\frac{1}{n}\sum_{i=1}^{n}\left| \frac{\pi_i^w\wedge 1}{\breve{\pi}_i^w\wedge 1}-1\right|\left[\hat y^{s}(\bm x_{i},\hat{\bm\theta})-y_i^{p}\right]^{2} \frac{\partial \hat y^{s}(\bm x_{i},\hat{\bm\theta})}{\partial\bm\theta} \frac{\partial \hat y^{s}(\bm x_{i},\hat{\bm\theta})}{\partial\bm\theta^T} \\
		\leq& \frac{\sum_{i=1}^{n}h_i^{mV}}{\sum_{i=1}^{n}h_0^{mV}(\bm x_i)} \sum_{i=1}^{n}\frac{\left[\hat y^{s}(\bm x_{i},\hat{\bm\theta})-y_i^{p}\right]^{2} \frac{\partial \hat y^{s}(\bm x_{i},\hat{\bm\theta})}{\partial\bm\theta} \frac{\partial \hat y^{s}(\bm x_{i},\hat{\bm\theta})}{\partial\bm\theta^T}}{n} \frac{\left|h_0^{mV}(\bm x_i)-h_i^{mV} \right|}{\rho r\frac{1}{n}\sum_{i=1}^{n}h_i^{mV}}  \\
		&+ \left|\frac{\sum_{i=1}^{n}h_i^{mV}}{\sum_{i=1}^{n}h_0^{mV}(\bm x_i)}-1 \right|  \sum_{i=1}^{n}\frac{\left[\hat y^{s}(\bm x_{i},\hat{\bm\theta})-y_i^{p}\right]^{2} \frac{\partial \hat y^{s}(\bm x_{i},\hat{\bm\theta})}{\partial\bm\theta} \frac{\partial \hat y^{s}(\bm x_{i},\hat{\bm\theta})}{\partial\bm\theta^T}}{n}  \frac{h_i^{mV}}{\rho r\frac{1}{n}\sum_{i=1}^{n}h_i^{mV}}.
	\end{aligned}
\end{equation}

Write $\frac{\partial \hat y^{s}(\bm x_{i},\hat{\bm\theta})}{\partial\bm\theta} \frac{\partial \hat y^{s}(\bm x_{i},\hat{\bm\theta})}{\partial\bm\theta^T}=\|\frac{\partial \hat y^{s}(\bm x_{i},\hat{\bm\theta})}{\partial\bm\theta^T}\|_E$. Now we bound (\ref{disV1}). 
By combining the triangle inequality and the consistency $\tilde{\bm\theta}_{0}$, it follows that
    \begin{equation*}
	\begin{aligned}
		&\left|h_0^{mV}(\bm x_i)-h_i^{mV} \right| \\
		\leq& \left|\hat y^{s}(\bm x_{i},\tilde{\bm\theta}_{0})-\hat y^{s}(\bm x_{i},\hat{\bm\theta}) \right| \left[ \frac{\partial\hat y^{s}(\bm x_{i},\tilde{\bm\theta}_{0})}{\partial\bm\theta}\tilde{\mathbf J}_0^{-2}\frac{\partial\hat y^{s}(\bm x_{i},\tilde{\bm\theta}_{0})}{\partial\bm\theta^{T}} \right]^{\frac{1}{2}} \\
		&+\left|\hat y^{s}(\bm x_{i},\tilde{\bm\theta}_{0})-\hat y^{s}(\bm x_{i},\hat{\bm\theta}) \right| \left[ \frac{\partial \hat y^{s}(\bm x_{i},\hat{\bm\theta})}{\partial\bm\theta} \tilde{\mathbf J}^{-2}\frac{\partial \hat y^{s}(\bm x_{i},\hat{\bm\theta})}{\partial\bm\theta^{T}} \right]^{\frac{1}{2}}.    
	\end{aligned}
\end{equation*}
 Assumption (H.\ref{a6}) and the consistency of $\tilde{\bm\theta}_0$ indicate that
\begin{align*}
	&\left| \hat{y}^{s}(\bm x_{i},\hat{\bm\theta})-\hat{y}^{s}(\bm x_{i},\tilde{\bm\theta}_0) \right| \leq m_{1}(x_i)\left\|\hat{\bm\theta}-\tilde{\bm\theta}_0 \right\|_E =o_p(1). 
\end{align*}
Also, by the consistency of $\tilde{\bm\theta}_0$, we have that
$\tilde{\mathbf{J}}_0$ converges to $\tilde{\mathbf J}$ in probability. 
Together with the fact that $ \left[ \frac{\partial \hat y^{s}(\bm x_{i},\hat{\bm\theta})}{\partial\bm\theta} \tilde{\mathbf J}^{-2}\frac{\partial \hat y^{s}(\bm x_{i},\hat{\bm\theta})}{\partial\bm\theta^{T}} \right]\leq \lambda_{\max}(\tilde{\mathbf J}^{-1})\left\|\frac{\partial\hat  y^{s}(\bm x_{i},\hat{\bm\theta})}{\partial\bm\theta^{T}} \right\|_E$, there is
\begin{align}
\label{disth}
\left|h_0^{mV}(\bm x_i)-h_i^{mV} \right|=o_p(1).
\end{align}
Thus,  there is
\begin{align}\label{Vdis1}
	&\frac{1}{n}\sum_{i=1}^{n}h_i^{mV}=\frac{1}{n}\sum_{i=1}^{n}h_0^{mV}(\bm x_i)+o_p(1).
\end{align}

Because 
\begin{equation*}
	\begin{aligned}
		\frac{1}{n}\sum_{i=1}^{n}h_i^{mV}&=\frac{1}{n}\sum_{i=1}^{n}\left| y_{i}^{p}-\hat{y}^{s}(\bm x_{i},\hat{\bm\theta})\right| \left[ \frac{\partial \hat{y}^{s}(\bm x_{i},\hat{\bm\theta})}{\partial\bm\theta}\tilde{\mathbf J}^{-2}\frac{\partial \hat{y}^{s}(\bm x_{i},\hat{\bm\theta})}{\partial\bm\theta^{T}} \right]^{1/2}\\
		&\geq \lambda_{\min}(\tilde{\mathbf J}^{-1})\frac{1}{n}\sum_{i=1}^{n}\left| y_{i}^{p}-\hat{y}^{s}(\bm x_{i},\hat{\bm\theta})\right| \left\| \frac{\partial \hat{y}^{s}(\bm x_{i},\hat{\bm\theta})}{\partial\bm\theta^{T}}\right\|_E,
	\end{aligned}
\end{equation*}
and 
$\frac{1}{n}\sum_{i=1}^{n}\left[ y_{i}^{p}-\hat{y}^{s}(\bm x_{i},\hat{\bm\theta})\right]^2 \left\| \frac{\partial \hat{y}^{s}(\bm x_{i},\hat{\bm\theta})}{\partial\bm\theta^{T}}\right\|_E^2=O_p(1)$ from Assumptions (H.\ref{a2}) – (H.\ref{a3}), it follows that
\begin{align}\label{Vdis2}
	&\left(\frac{1}{n}\sum_{i=1}^{n}h_i^{mV}\right)^{-1}=O_p(1).
\end{align}

By plugging (\ref{disth}) – (\ref{Vdis2}) into (\ref{disV1}), there is
\begin{equation}
	\left\| \breve{\mathbf V}^w-\tilde{\mathbf V}^w \right\|_s=o_p(r^{-1}). \label{dl5}
\end{equation}

\begin{itemize}
	\item Next, we consider $\frac{\partial^2 l_{\tilde{\bm\theta}_{0}}^{*}(\hat{\bm\theta})}{\partial\bm\theta^T \partial \theta_j}$,	which can be written as
\end{itemize}
\begin{equation} 
	\label{ddl2}
	\begin{aligned}	
		\frac{\partial^2 l_{\tilde{\bm\theta}_{0}}^{*}(\hat{\bm\theta})}{\partial\bm\theta^T \partial \theta_j}
		=\frac{1}{n}\frac{\partial^{2}\sum_{i=1}^{n}\frac{a_{i}}{\breve{\pi}_i^w\wedge 1}\left[y^p_i-\hat{y}^{s}(\bm x_{i},\hat{\bm\theta})\right]^{2}}{\partial\bm\theta^T \partial\theta_j}.
	\end{aligned}
\end{equation}		
Using the arguments similar to (\ref{515}), (\ref{ddl2}) can be represented by
\begin{equation}
	\frac{1}{n}\frac{\partial^{2}\sum_{i=1}^{n}\frac{a_{i}}{\breve{\pi}_i^w\wedge 1}\left[y^p_i-\hat{y}^{s}(\bm x_{i},\hat{\bm\theta})\right]^{2}}{\partial\bm\theta^T \partial\theta_j}=\frac{1}{n}\sum_{i=1}^{n} \frac{\partial^{2}\left[y^p_i-\hat{y}^{s}(\bm x_i,\hat{\bm\theta})\right]^{2}}{\partial\bm\theta^T \partial\theta_j} + O_{p}(r^{-\frac{1}{2}}).
	\label{ddll}
\end{equation}

\begin{itemize}
	\item Third, we bound $\int_{0}^{1}\int_{0}^{1}  
	\frac{\partial^3 l_{\tilde{\bm\theta}_{0}}^{*}(\hat{\bm\theta}+(\breve{\bm\theta}-\hat{\bm\theta})uv)}{\partial\bm\theta \partial\bm\theta^T \partial\theta_j}vdudv$.
\end{itemize}
Suppose $\bm\theta_0\in \bm\Theta_0$, next we bound $\frac{\partial^3 l_{\tilde{\bm\theta}_{0}}^{*}(\bm\theta_0)}{\partial\theta_k \partial\theta_{k'} \partial\theta_j}$, which can be written as follows
\begin{equation*} 
	\begin{aligned}	
		\frac{\partial^3 l_{\tilde{\bm\theta}_{0}}^{*}(\bm\theta_0)}{\partial\theta_k \partial\theta_{k'} \partial\theta_j}
		&=\frac{1}{n}\frac{\partial^{3}\sum_{i=1}^{n}\frac{a_{i}}{\breve{\pi}_i^w\wedge 1}\left[y^p_i-\hat y^{s}(\bm x_{i},{\bm\theta}_0)\right]^{2}}{\partial\theta_k\partial \theta_{k'} \partial \theta_j}.
	\end{aligned}
\end{equation*}	
Similar with (\ref{515}), $\frac{\partial^3 l_{\tilde{\bm\theta}_{0}}^{*}(\bm\theta_0)}{\partial\theta_k \partial\theta_{k'} \partial\theta_j}$ can be bounded by
 $$\frac{1}{n}\frac{\partial^{3}\sum_{i=1}^{n}\left[y^p_i-\hat y^{s}(\bm x_{i},{\bm\theta}_0)\right]^{2}}{\partial\theta_k\partial \theta_{k'} \partial \theta_j}+O_p(r^{-1/2}).$$
 Because $\|\breve{\bm\theta}-\hat{\bm\theta}\|_E=o_{p|\mathcal D_n}(1)$, $\hat{\bm\theta}+(\breve{\bm\theta}-\hat{\bm\theta})uv\in \bm\Theta_0$, similar to (\ref{dddl3}), it follows that
\begin{equation}
	\begin{aligned}
		\int_{0}^{1}\int_{0}^{1}  
		\frac{\partial^3 l_{\tilde{\bm\theta}_{0}}^{*}(\hat{\bm\theta}+(\breve{\bm\theta}-\hat{\bm\theta})uv)}{\partial\theta_k \partial\theta_k^{'} \partial\theta_j}vdudv 
		=&  \int_{0}^{1}\int_{0}^{1}  \frac{1}{n}\frac{\partial^{3}\sum_{i=1}^{n}\left[y^p_i-\hat y^{s}(\bm x_{i},{\bm\theta}_0)\right]^{2}}{\partial\theta_k\partial \theta_{k'} \partial \theta_j}
		  vdudv  + O_p(r^{-\frac{1}{2}})\\
		=&O_p(1).
		  \label{dddl2}
	\end{aligned}
\end{equation}

The desired result can be obtained by plugging (\ref{dl2}) and (\ref{dl5}) – (\ref{dddl2}) into (\ref{taylor2}).



\end{appendices}

\bibliographystyle{chicago}

\bibliography{reference}

\end{document}